  \documentclass{amsart}
  \usepackage{url}
  \usepackage{amssymb}
  
  \input xypic
  \xyoption{all}
  
  \usepackage{amsmath, amsthm,  amsfonts}

  \theoremstyle{plain}
  \newtheorem{theorem}{Theorem}[section]
  \newtheorem{lemma}[theorem]{Lemma}
  \newtheorem{proposition}[theorem]{Proposition}

  \theoremstyle{definition}
  \newtheorem{definition}[theorem]{Definition}

  \theoremstyle{remark}
  \newtheorem{remark}[theorem]{Remark}
  
  \newtheorem{note}{Note}[section]

  \numberwithin{equation}{section}
  \numberwithin{figure}{section}

  \newcommand{\cQ}{{\mathcal Q}}
  \newcommand{\cA}{{\mathcal A}}

  \renewcommand{\cL}{{\mathcal L}}
  \renewcommand{\cD}{{\mathcal D}}
  \newcommand{\cP}{{\mathcal P}}
  \newcommand{\cU}{{\mathcal U}}
  \newcommand{\cG}{{\mathcal G }}

   \newcommand{\sta}{\stackrel}
   \newcommand{\ba}{\begin{eqnarray}}
   \newcommand{\na}{\end{eqnarray}}
   \newcommand{\ban}{\begin{eqnarray*}}
   \newcommand{\nan}{\end{eqnarray*}}

  \newcommand{\g}{{\mathfrak g}}

  \newcommand{\CC}{{\mathbb C}}
  \newcommand{\RR}{{\mathbb R}}
  \newcommand{\ZZ}{{\mathbb Z}}
  \newcommand{\QQ}{{\mathbb Q}}
  \renewcommand{\AA}{{\mathbb A}}

  \renewcommand{\a}{\alpha}
  \renewcommand{\b}{\beta}
  \renewcommand{\c}{\gamma}
  
  \renewcommand{\d}{\delta}

    \newcommand{\disp}{\displaystyle}

\begin{document}

  \title[Bundle gerbes for Chern-Simons and Wess-Zumino-Witten theories]
{Bundle gerbes for Chern-Simons and Wess-Zumino-Witten theories}

  \author[A.L. Carey]{Alan L. Carey}
\address[Alan L. Carey]
  {Mathematical Sciences Institute\\
  Australian National University\\
  Canberra ACT 0200 \\
  Australia}
  \email{acarey@maths.anu.edu.au}

  \author[S. Johnson]{Stuart Johnson}
  \address[Stuart Johnson]
  {Department of Pure Mathematics\\
  University of Adelaide\\
  Adelaide, SA 5005 \\
  Australia}
  \email{sjohnson@maths.adelaide.edu.au}

  \author[M.K. Murray]{Michael K. Murray}
  \address[Michael K. Murray]
  {Department of Pure Mathematics\\
  University of Adelaide\\
  Adelaide, SA 5005 \\
  Australia}
  \email{mmurray@maths.adelaide.edu.au}

  \author[D. Stevenson]{Danny Stevenson}
  \address[Danny Stevenson]
  {Department of Mathematics\\
202 Surge Building\\
University of California, Riverside\\
Riverside CA 92521-0135\\
USA
}
  \email{dstevens@math.ucr.edu}

  \author[Bai-Ling Wang]{Bai-Ling Wang}
  \address[Bai-Ling Wang]
  {Institut f\"ur Mathematik\\
    Universit\"at Z\"urich\\
   Winterthurerstrasse 190\\
   CH-8057, Z\"urich\\
 Switzerland}
  \email{bwang@math.unizh.ch}

  \thanks{The authors acknowledge the support of the Australian
  Research Council. ALC thanks  MPI f\"ur Mathematik in Bonn and ESI in Vienna and 
  BLW thanks CMA  of Australian National University for
  their hospitality during part of the writing of this paper. }

  \begin{abstract} We develop the theory of Chern-Simons bundle 2-gerbes and
   multiplicative bundle gerbes associated to any
principal $G$-bundle with connection
   and a class in $H^4(BG, \ZZ)$ for  a compact semi-simple
Lie group $G$. The Chern-Simons bundle 2-gerbe realises differential geometrically
 the Cheeger-Simons invariant.  We apply these notions to refine
   the  Dijkgraaf-Witten correspondence
   between three dimensional Chern-Simons functionals
and Wess-Zumino-Witten models
   associated to the group $G$. We do this by introducing a
 lifting to the level of bundle gerbes of the natural map from
   $H^4(BG, \ZZ)$ to $H^3(G, \ZZ)$.
   The notion of a multiplicative bundle gerbe
   accounts geometrically for the subtleties in this correspondence
   for non-simply connected  Lie groups.
The implications for Wess-Zumino-Witten models
   are also discussed.

  \end{abstract}
  \maketitle

\tableofcontents

  \section{Introduction}

In \cite{Qui} Quillen introduced the determinant line bundle of
Cauchy-Riemann operators on a Hermitian vector bundle coupled to
unitary connections over a Riemann surface. This work influenced
the development of many lines of investigation including the study
of Wess-Zumino-Witten actions on Riemann surfaces.  Note that
Quillen's determinant line bundle also plays an essential role in
the construction of the universal bundle gerbe in \cite{CW}, see
also \cite{Bry2}.

The relevance of Chern-Simons gauge theory has been noted
by many authors, starting with Ramadas-Singer-Weitsman \cite{RSW} and
recently Dupont-Johansen \cite{DupJoh},
who used gauge covariance of the Chern-Simons functional to give a
geometric construction of Quillen line bundles. The curvatures of these
line bundles in an analytical set-up were studied extensively by
Bismut-Freed \cite{BisFre} and in dimension two, went back to the
Atiyah-Bott work on the Yang-Mills equations over Riemann surfaces.
\cite{AB}.

A new element was introduced into this picture by Freed \cite{Fre0} and
\cite{Fre} (a related
line of thinking was started by some of the present authors \cite{CMW})
through the introduction of higher algebraic structures
(2-categories) to study Chern-Simons functionals on
3-manifolds with boundary and corners. For  closed 3-manifolds
one needs to study the behaviour of the Chern-Simons action under
gluing formulae (that is topological quantum field theories)
generalising the corresponding picture for Wess-Zumino-Witten.
Heuristically, there is
a Chern-Simons line bundle as in \cite{RSW}, such that for a 3-manifold
with boundary,
the Chern-Simons action is a section of the Chern-Simons line bundle
associated to the boundary Riemann
surface.  For a codimension two submanifold, a closed circle, the
Chern-Simons action takes values in a $U(1)$-gerbe or an abelian
group-like 2-category.

Gerbes  first began to enter the picture with J-L Brylinski \cite{Bry0}
and Breen \cite{Bre}. The latter
 developed the notion of a 2-gerbe as a sheaf of bicategories
extending Giraud's \cite{Gir} definition of a gerbe as a sheaf of groupoids.
J-L Brylinski used Giraud's gerbes to study the central extensions of loop
groups, string structures and the relation to Deligne cohomology. With
McLaughlin, Brylinski developed a  2-gerbe over a manifold $M$ to realise
degree 4 integral cohomology on $M$ in \cite{BryMcL1}
and introduced an  expression of the 2-gerbe holonomy as a
Cheeger-Simons differential character on any manifold with a
triangulation. This is the starting point for Gomi \cite{Gom1},\cite{Gom2}
who developed a local theory of the Chern-Simons functional
along the lines of Freed's suggestion.
A different approach to some of
these matters using simplicial manifolds has been found by Dupont and Kamber
\cite{DupKam}.

Our contribution is to develop a global differential geometric
realization of Chern-Simons functionals
using a Chern-Simons bundle
2-gerbe and to apply this to the question raised by Dijkgraaf and Witten
about the relation between Chern-Simons and Wess-Zumino-Witten models.
Our approach provides a unifying perspective on
 all of this previous work in a fashion
that can be directly related to the physics literature on Chern-Simons
field theory (thought of as a path integral defined in terms of the
Chern-Simons functional).

  In \cite{DijWit} it is shown that three dimensional Chern-Simons gauge
   theories with gauge group $G$ can be classified by the integer
cohomology group
   $H^4(BG, \ZZ)$, and
  conformally invariant sigma models in two dimension with target space a compact Lie group
  (Wess-Zumino-Witten models)
can be classified by $H^3(G, \ZZ)$. It  is  also
 established that  the correspondence between three dimensional Chern-Simons
gauge theories and
 Wess-Zumino-Witten models   is related to the transgression map
    $$\tau: 
  H^4(BG, \ZZ) \to H^3(G, \ZZ),
  $$
  which explains the subtleties in this correspondence for compact,
  semi-simple  non-simply connected  Lie groups (\cite{MooSeib}).

  In the present work we introduce Chern-Simons bundle 2-gerbes
  and the notion of  multiplicative bundle gerbes,  and
  apply them to explore the geometry of the Dijkgraaf-Witten
  correspondence. To this end,
we will assume throughout that $G$ is a compact semi-simple Lie group.

 The role of Deligne cohomology as an ingredient in
topological field theories goes back to \cite{Gaw}
and we add a new feature in section \ref{Deligne} by  using
   Deligne cohomology valued characteristic classes for principal
  $G$-bundles with connection. Briefly speaking, a
  degree $p$ Deligne characteristic class for principal $G$-bundles with
connection is
  an  assignment to any principal $G$-bundle  with connection
    over $M$ of a  class in the degree $p$ Deligne cohomology group
$H^p(M, \cD^p)$  satisfying a certain functorial property.
 Deligne cohomology valued characteristic classes refine
 the characteristic classes for principal $G$-bundles.

  We will define three dimensional Chern-Simons gauge theories $CS(G)$
  as degree 3 Deligne cohomology valued characteristic classes for principal
  $G$-bundles with connection, but will later show that there is a global
differential geometric structure, the Chern-Simons bundle 2-gerbe,
associated to each Chern-Simons gauge theory.
 We will interpret a
Wess-Zumino-Witten model as arising from the curving of
a bundle gerbe associated to a
degree $2$ Deligne cohomology class on the Lie group $G$
as in \cite{CMM} and
\cite{GawRei}.
   We then use a
  certain canonical $G$ bundle defined on $S^1 \times G$ to
  construct a transgression map between classical Chern-Simons
gauge theories $CS(G)$
   and classical  Wess-Zumino-Witten
 models $WZW(G)$ in section \ref{section:CS2WZW},
which is a lift of the transgression map
    $H^4(BG, \ZZ) \to H^3(G, \ZZ)$. The resulting correspondence
    $$
    \Psi:  CS(G) \longrightarrow WZW(G)
    $$
    refines the  Dijkgraaf-Witten correspondence between three dimensional Chern-Simons gauge
   theories and Wess-Zumino-Witten models associated to a compact Lie group
   $G$. On Deligne cohomology groups,
our correspondence $\Psi$ induces a transgression
   map
   $$
   H^3(BG, \cD^3) \longrightarrow H^2(G, \cD^2),
   $$ and  refines the natural transgression map $\tau: H^4(BG, \ZZ) \to H^3(G, \ZZ)$
   (Cf. Proposition \ref{true:trans}).
   See \cite{BryMcL} for a related transgression of Deligne cohomology in a different set-up.

   For any integral cohomology class in $H^3(G, \ZZ)$, there is a unique stable
 equivalence class of bundle gerbe (\cite{Mur, MurSte0}) whose Dixmier-Douady class
 is the given degree 3 integral cohomology class.
 Geometrically $H^4(BG, \ZZ)$ can be regarded as stable equivalence
  classes of bundle $2$-gerbes over $BG$, whose induced bundle gerbe over
$G$ has a certain  multiplicative structure.

   To study the geometry of the correspondence $\Psi$, we revisit the
   bundle 2-gerbe theory developed in \cite{Ste} and \cite{Joh} in Section \ref{b2g}.
   Note that transformations between stable isomorphisms provide $2$-morphisms making 
the category $\mathbf{BGrb}_M$ of bundle gerbes over $M$
and stable isomorphisms between bundle gerbes into a  
\emph{bi-category} (Cf. \cite{Ste}). 

For a smooth surjective submersion $\pi\colon X\to M$, 
consider the face operators $\pi_i: X^{[n]} \to X^{[n-1]}$ on the simplicial 
manifold $X_\bullet = \{X_n = X^{[n+1]}\}$. 
Then a bundle $2$-gerbe on $M$ consists of the data of a 
smooth surjective submersion $\pi\colon X\to M$ together 
with
\begin{enumerate}
\item  An object $(\cQ,Y,X^{[2]})$
in $\mathbf{BGrb}_{X^{[2]}}$.
\item  A stable isomorphism $m$: $\pi_1^*\cQ\otimes \pi_3^*\cQ  
\to \pi_2^*\cQ$ in $\mathbf{BGrb}_{X^{[3]}}$ defining the bundle 2-gerbe
product which is associative up to a 2-morphism $\phi$
in $\mathbf{BGrb}_{X^{[4]}}$. 
\item The 2-morphism $\phi$ satisfies a natural coherency
condition in $\mathbf{BGrb}_{X^{[5]}}$.
\end{enumerate}
   
   We then develop a multiplicative bundle gerbe theory
   over $G$ in section \ref{section:multiplicative-gerbe} as a
   simplicial bundle gerbe on the simplicial manifold associated to $BG$.
   We say a bundle gerbe $\cG$ over $G$ is {\it transgressive} if   
   the Deligne class of $\cG$, written  $d(\cG)$ is in the image of
  the correspondence map $\Psi: CS(G) \to WZW(G) = H^2(G, \cG^2)$.
   The main results of this paper are
 the following two theorems  (Theorem \ref{trans=multi:chern}
   and Theorem \ref{CS:image})

   \begin{enumerate}
   \item {\sl  The Dixmier-Douady class of a bundle gerbe $\cG$ over $G$  lies in 
   the image of  the transgression map
  $\tau: H^4(BG,\ZZ) \rightarrow H^3(G,\ZZ)$ if and only if
$\cG$ is multiplicative.}
  \item {\sl  Let $\cG$ be a bundle
    gerbe over $G$ with connection and curving,
whose Deligne class $d(\cG)$ is in $H^2(G, \cD^2)$.
 Then $\cG$ is transgressive if and only if $\cG$ is multiplicative. }
    \end{enumerate}

   Let $\phi$ be an element in $H^4(BG, \ZZ)$. The corresponding $G$-invariant
   polynomial on the Lie algebra under the universal Chern-Weil homomorphism is denoted
   by $\Phi$. For any connection $\AA$ on the universal bundle $EG\to BG$
   with the curvature form $F_\AA$,
   \[
   (\phi, \Phi(\disp{\frac{i}{2\pi}}F_\AA)) \in H^4(BG, \ZZ) \times_{H^4(BG, \RR)}
   \Omega^4_{cl, 0}(BG)
   \]
   (where $\Omega^4_{cl, 0}(BG)$ is the space of closed 4-forms on $BG$ with
   periods in $\ZZ$),
defines a unique  degree 3 Deligne class in  $H^3(BG, \cD^3)$.  Here
   we fix a smooth infinite dimensional model of $EG\to BG$
   by embedding $G$ into $U(N)$ and letting $EG$ be the Stiefel manifold of
    $N$ orthonormal vectors in a separable complex Hilbert space.

   We will show that $H^3(BG, \cD^3)$ classifies the stable
equivalence classes of
   bundle 2-gerbes with curving on $BG$, (we already know that the second
   Deligne cohomology classifies the stable  equivalence classes of
  bundle gerbes with curving).  These are the
   universal Chern-Simons bundle 2-gerbes $\cQ_\phi$ in section \ref{CS2gerbe}
   (cf. Proposition \ref{cs2bg:exist}) giving a geometric realisation of
   the degree 3 Deligne class
   determined by $(\phi, \Phi(\disp{\frac{i}{2\pi}}F_\AA))$.

   We show that
 for any principal $G$-bundle $P$ with connection $A$ over $M$,
the associated Chern-Simons bundle
   2-gerbe $\cQ_\phi(P, A)$ over $M$ is obtained by the pull-back
of the universal
    Chern-Simons bundle 2-gerbe $\cQ_\phi$
   via a classifying map. {\sl The bundle 2-gerbe curvature of $\cQ_\phi(P, A)$
   is given by $\Phi(\disp{\frac{i}{2\pi}}F_A)$, and the bundle 2-gerbe curving
   is given by the Chern-Simons form associated to $(P, A)$ and $\phi$}.

  Under the canonical isomorphism between Deligne cohomology and
   Cheeger-Simons cohomology,
there is a canonical holonomy map for any degree
  $p$ Deligne class from the group of smooth $p$-cocycles to $U(1)$.
   This holonomy is known as the
   Cheeger-Simons differential character associated to the Deligne class.

The bundle 2-gerbe holonomy for this Chern-Simons bundle
    2-gerbe $\cQ_\phi(P, A)$ over $M$ as given by the Cheeger-Simons
    differential character is used in the integrand for
 the path integral for the Chern-Simons
quantum field theory. In the $SU(N)$ Chern-Simons theory, $\Phi$ is chosen to
  be the second Chern polynomial.  For a smooth map $\sigma: Y \to M$,
  under a fixed trivialisation of
  $\sigma^*(P, A)$ over $Y$,  the corresponding  holonomy of $\sigma$   is given by $e^{2\pi i CS(\sigma, A)}$, where $CS(\sigma, A)$ can be
  written as  the following well-known Chern-Simons form:
  \[
   \disp{\frac{k}{8\pi^2}\int_Y} Tr\sigma^*(A\wedge dA + \disp{\frac 13} A\wedge A \wedge A),
 \]
   Here $k\in \ZZ$ is the level determined by
$\phi \in H^4(BSU(N), \ZZ) \cong \ZZ$.

We will establish in Theorem \ref{csb2g=cs-invariant} that
{\sl the Chern-Simons bundle  2-gerbe $\cQ_\phi(P, A)$ over $M$ is equivalent in 
Deligne cohomology to the  Cheeger-Simons invariant associated to the principal $G$-bundle $P$ with
   a connection $A$ and   a class $\phi \in H^4(BG, \ZZ)$.}

  In the concluding section we
 show that the Wess-Zumino-Witten models in the
  image of this correspondence $\Psi$ satisfy a quite interesting
 multiplicative property that is associated  with the group multiplication
on $G$. This multiplicative property is a feature of the holonomy of
every multiplicative bundle gerbe.  It implies that for the
transgressive Wess-Zumino-Witten models, the
  so-called $B$ field satisfies a certain integrality condition.
   Using our multiplicative bundle gerbe theory,  we
can give a very satisfying
explanation why, for non-simply connected groups, multiplicative
  bundle gerbes only exist for Dixmier Douady classes that
are certain particular multiples of the generator
  in  $H^3(G, \ZZ)$ (we call this the `level'). For a non-simply connected Lie
  group, there exists a subtlety in the construction of positive energy
  representations of its loop group, see \cite{PS} \cite{Tol}, where the level 
  is defined in terms of its Lie algebra. 
  
  While this paper was in preparation,  Aschieri and Jur\v{c}o  in \cite{AscJur}
 proposed a similar construction of Chern-Simons 2-gerbes in terms of
 Deligne classes developed in \cite{Joh} to study M5-brane anomalies and 
 $E_8$ gauge theory. Their discussions have some overlaps with our local 
 descriptions of bundle 2-gerbes and holonomy of 2-gerbes. 

  \section{Deligne characteristic classes for principal $G$-bundles}\label{Deligne}

 In this Section, we first review briefly  Deligne and Cheeger-Simons
 cohomology, and then define a Deligne cohomology valued characteristic class for any
 principal $G$-bundle with connection  over a smooth manifold $M$
 with $G$ a compact semi-simple Lie group.

  Let $H^p(M, \cD^p)$ be the $p$-th Deligne cohomology group,
   which is the hypercohomology group of the complex of sheaves on $M$:
  \[
  \underline{U(1)}  \stackrel{d \text{log}} {\rightarrow} \Omega^1_M
  \stackrel{d} {\rightarrow} \cdots \stackrel{d} {\rightarrow}  \Omega^p_M
  \]
  where $\underline{U(1)}$ is the sheaf of smooth $U(1)$-valued functions on $M$,
  $\Omega^p_M$  is the sheaf of imaginary-valued differential $p$-forms on $M$.
    Take any degree $p$ Deligne class
    \[
    \xi = [ g, \omega^1, \cdots \omega^p],
    \]
 then with respect to a good cover of $M$, $\{g_{i_0i_1\cdots i_p}\}$
  represents a $U(1)$-valued \v{C}ech
    $p$-cocycle on $M$, and hence defines an element in
    \[
    H^p(M, \underline{U(1)}) \cong H^{p+1}(M, \ZZ).
    \]
 The corresponding element in $H^{p+1}(M, \ZZ)$ is denoted by $c(\xi)$,
and referred to as {\sl the characteristic  class} of $\xi$. Moreover,
    $d \omega^p$ is a globally defined closed $p+1$ form on $M$ with periods
    in $ 2\pi i  \ZZ$  called {\sl the
    curvature} of $\xi$ and denoted by $curv(\xi)$. Without causing any confusion, we often
    identify $curv(\xi)$ with $curv(\xi)/2\pi i$ whose
     periods are in $\ZZ$.

  We   have indexed the Deligne cohomology group so that a degree $p$
Deligne class
  has holonomy (to be discussed later in this section) over $p$ dimensional sub-manifolds and a characteristic
  class in $H^{p+1}(M, \ZZ)$. For example, $H^1(M, \cD^1)$ is the space
  of  equivalence
  classes  of line bundles with connection, whose holonomy is defined for any
  smooth path and whose characteristic
  class in $H^{2}(M, \ZZ)$ is given by the first Chern class of the underlying line
  bundle. Next in the hierarchy, $H^2(M, \cD^2)$ is the space of stable
  isomorphism classes of bundle gerbes with connection and curving,
whose holonomy
  is defined for any 2-dimensional closed sub-manifold and whose characteristic
  class in $H^{3}(M, \ZZ)$ is given by the Dixmier-Douady class of the underlying bundle
  gerbe.

  The  Deligne cohomology group $H^p(M, \cD^p)$ is part of
 the following  two  exact sequences (Cf. \cite{Bry})
  \ba\label{exact:1}
  0 \rightarrow \Omega^p_{cl, 0}(M) \rightarrow
  \Omega^p(M) \rightarrow H^p(M, \cD^p) \stackrel{c}{\rightarrow} H^{p+1}(M, \ZZ)
 \rightarrow 0
  \na
  where $\Omega^p_{cl, 0}(M) $ is the subspace of closed
  $p$-forms on $M$ with periods in $\ZZ$, in the space of
  $p$-forms $\Omega^p(M)$, and $c$ is the  characteristic class  map; and
  \ba\label{exact:2}
  0 \rightarrow H^{p}(M, \RR/\ZZ) \rightarrow
  H^p(M, \cD^p) \stackrel{curv}{\rightarrow}  \Omega^{p+1}_{cl, 0}(M)
  \rightarrow 0
  \na
  where the map $curv$ is the curvature map on $H^p(M, \cD^p)$.

  We remark  that
  for a Deligne class $\xi \in H^p(M, \cD^p)$, $curv(\xi)$ the curvature of $\xi$
   defines   a  cohomology class in $H^{p+1}(M, \RR)$
    which agrees with the image of
  $c(\xi)$ under the map $H^{p+1}(M, \ZZ) \to H^{p+1}(M,
  \RR)$ sending an integral class to a real class.

  Recall that the Cheeger-Simons group of differential characters of degree
  $p$ in  \cite{CheSim}, $\check{H}^p(M, U(1))$, is defined to be the space of pairs, $(\chi, \omega)$
  consisting of a homomorphism
  \[
  \chi:  Z_p(M, \ZZ) \to U(1)
  \]
  where $Z_p(M, \ZZ)$ is the group of   smooth $p$-cycles, and
  an  imaginary-valued closed $(p+1)$-form $\omega$ on $M$ with periods in
  $2\pi i \ZZ$ such that  for any smooth $(p+1)$ chain $\sigma$
  \[
  \chi(\partial \sigma) = exp(\disp{\int_\sigma} \omega).
  \]
  The Cheeger-Simons  group $\check{H}^p(M, U(1))$ enjoys the
  same exact sequences (\ref{exact:1}) and (\ref{exact:2}) as the Deligne cohomology
  group $H^p(M, \cD^p)$. In fact, the holonomy and the curvature of a Deligne class
  define  a canonical isomorphism
  \ba\label{del=cs}
  (hol, curv):  H^p(M, \cD^p) \longrightarrow \check{H}^p(M, U(1)).
  \na
  Here the holonomy of a Deligne class $\xi$ is defined as follows.
 For a smooth
  $p$-cycle given by a triangulation of
  a smooth map $X \to M$, pull back $\xi$ to $X$ to obtain a Deligne
  class  on $X$. Lift this class to an element $\alpha$ in $\Omega^p(X)$ from
  the exact sequence (\ref{exact:1}) as $H^{p+1}(X, \ZZ) =0$ and then
  \[
   hol(\xi) =  exp(\disp{\int_X} \a)
  \]
  is independent of the choice of $\a$, this again follows
  from (\ref{exact:1}). For a general smooth $p$-cycle $\sigma
  =\sum_k n_k\sigma_k$, we choose
  a local representative $(g, \omega^1, \cdots \omega^p)$ of $\xi$
  under a good cover $\{U_i\}$ of $M$ such that for each  smooth $p$-simplex $f_k:
  \sigma_k\to M$,
  with possible subdivisions, $f_k(\sigma_k)$
    is contained in some open set $U_{i_k}$ for which
    $f_k$ has a smooth extension.
 We can define (see \cite{Gaw} for $p=2$, and \cite{CMJ, Gom1} for $p=2,3$)
   \[\begin{array}{lll}
  \cA(\sigma_k) &=& exp\Bigl\{\disp{\int_{\sigma_k}f^*_k\omega^p_{i_k}
  + \sum_{\tau_{(1)}}\int_{\tau_{(1)}\subset \sigma_k}f_k^* \omega^{p-1}_{i_{\tau_{(1)}} i_k}
  + \cdots  }\Bigr\}\\[2mm]
  && \cdot \prod_{\tau_{(p)}\in \tau_{(p-1)}\subset \cdots \subset \tau_{(1)}\subset
  \sigma_k}
    g^{\ }_{i_k i_{\tau_{(1)}} \cdots i_{\tau_{(p)}}}
    \end{array}
    \]
  where $\{\tau_{(j)}\}$ is the set of codimension $j$ faces of $\sigma_k$,
  with $\tau_{(j)}\subset U_{i_{\tau_{(j)}}}$ and induced orientation
  from $\sigma_k$. It is routine to show that
  \[
  hol(\xi) =\prod_k \cA (\sigma_k)^{n_k}
  \]
  is independent of local representative of $\xi$ and subdivisions because
  $\sum_k n_k\sigma_k$ is a cycle.

  With the understanding  of (\ref{del=cs}),
we often identify a Deligne class with the
  corresponding Cheeger-Simons differential character.  In \cite{HopSin},
  Hopkins and Singer develop a cochain model for the Cheeger-Simons
  cohomology where certain integrations are  well-defined,
  see also \cite{DupLju}. In this paper,
we need to define an integral on the Deligne cohomology group:
  \ba\label{integration:1}
  \disp{\int_{S^1}}: \qquad H^3(S^1\times G, \cD^3) \to H^2(G, \cD^2).
  \na
While this makes sense under the map $(hol, curv)$ and Hopkins-Singer's integration
  for the cochain model for the Cheeger-Simons differential characters
we will instead apply
  the exact sequences (\ref{exact:1}) and (\ref{exact:2}) to uniquely
   define  the integration map
  (\ref{integration:1})
  via the following  two commutative diagrams:
  \ba\label{integration:2}
  \begin{array}{ccccc}
  H^{3}(S^1\times G, \RR/\ZZ) &\to &
  H^{3}(S^1\times G, \cD^3) &\stackrel{curv}{\to} & \Omega^{4}_{cl, 0}(S^1\times G)\\[2mm]
  \downarrow\int_{S^1}&&\downarrow\int_{S^1}&&\downarrow\int_{S^1}\\[2mm]
  H^{2}(G, \RR/\ZZ) &\to &
  H^2(G, \cD^2) &\stackrel{curv}{\to} & \Omega^{3}_{cl, 0}(G),
  \end{array}
  \na
  where the integration map
  (\ref{integration:1}) is well-defined modulo the image of
  $H^{2}(G, \RR/\ZZ) \to H^2(G, \cD^2)$;
  and
    \ba
    \label{integration:3}
  \begin{array}{ccccc}
  \disp{\frac{\Omega^{3}(S^1\times G)}
  { \Omega^{3}_{cl, 0}(S^1\times G)}} &\to &
  H^{3}(S^1\times G, \cD^3) &\stackrel{c}{\to} & H^{4}(S^1\times G, \ZZ)\\[2mm]
  \downarrow\int_{S^1}&&\downarrow\int_{S^1}&&\downarrow\int_{S^1}\\[2mm]
  \disp{\frac{\Omega^{2}(G)}
  { \Omega^{2}_{cl, 0}(G)}}&\to &
  H^2(G, \cD^2) &\stackrel{c}{\to} & H^{3}(G, \ZZ),
  \end{array}
  \na
  where the integration map
  (\ref{integration:1}) is well-defined modulo the image of
  $\Omega^{2}(G)\to H^2(G, \cD^2)$.

To see how these two exact sequences may be used to uniquely specify the
integration map we let
$r \colon H^p(M, \ZZ) \to H^p(M,
  \RR)$ be the map
  that sends an integral class to a real class.
  Cheeger and Simons \cite{CheSim} define
  \ba\label{coincident}
  R^p(M, \ZZ) = \{(\omega, u) \in
  \Omega^p_{cl,0}(M) \oplus H^p(M, \ZZ) \mid r(u) = [\omega] \}
  \na
  where
  $[\omega]$ is the real cohomology class of the differential form
  $\omega$.
This  enables them to combine the map of a Deligne class to
  its characteristic class
  and the map to its curvature into one map from the Deligne cohomology
 group  into $R^{p+1}$.
  There is a short exact sequence:
  \ba\label{cs:coinc}
  0 \to
  \frac{H^p(M, \RR)}{r(H^p(M, \ZZ))} \to H^p(M, \cD^p)
  \sta{(curv, c)}{\longrightarrow} R^{p+1} (M, \ZZ) \to
  0.
  \na
 Then the following induced diagram is commutative:
  \[
  \begin{array}{cccccc}
  \disp{\frac{H^{3}(S^1\times G, \RR)}
  { r(H^{3}(S^1\times G, \ZZ))}} &\to &
  H^{3}(S^1\times G, \cD^3) &\stackrel{(curv, c)}{\longrightarrow} 
  & R^{4}(S^1\times G, \ZZ)&\\[2mm]
  \downarrow\int_{S^1}&&\downarrow\int_{S^1}&&\downarrow\int_{S^1}&\\[2mm]
  0 \quad &\to &
  H^2(G, \cD^2) &\stackrel{(curv, c)}{\longrightarrow} & R^{3}(G, \ZZ)&\to 0.
  \end{array}
  \]
  This commutative diagram and the fact that
  $H^2(G, \RR)=0$ for any compact
  semi-simple Lie group $G$ is the reason  the integration map
(\ref{integration:1}) is well defined.

  With these preparations we may now define Deligne characteristic classes
   for principal $G$-bundles with connection.
  Recall  that    a characteristic class $c$ for principal $G$-bundles is
  an assignment of a class $c(P)\in H^*(M, \ZZ)$ 
  to every isomorphism class of principal $G$-bundles
  $P \to M$.  Of course we could do
  $\QQ$, $\RR$ etc instead of $\ZZ$. This assignment is required to be
  `functorial' in the following sense: if $f \colon N \to M$ is a smooth
  map we require that $c(f^{*}(P)) = f^*(c(P))$, where $f^*P$ is the pull-back
  principal $G$-bundle over $N$.

  It is a standard fact that characteristic classes are in bijective correspondence with
  elements of $H^*(BG, \ZZ)$. The proof is: given 
  a  characteristic class $c$,  we  have, of course,
  $c(EG) \in H^*(BG, \ZZ)$ and conversely   if $\xi \in H^*(BG,
  \ZZ)$ is given, then defining $c_\xi(P) = f^*(\xi)$
  for any  classifying map $f \colon M \to BG$ gives rise to a characteristic class
  for the isomorphism class of principal $G$-bundles defined by the classifying map $f$.
  This uses the fact  that any two classifying maps $f$ and $g$ are
  homotopy equivalent so that $f^* = g^*$ and hence $f^*(\xi) = g^*(\xi)$.
  This definition motivates our definition of Deligne characteristic classes
  for principal $G$-bundles with connection.

  Note that the characteristic classes only depend
  on the underlying topological principal $G$-bundle, in order to define
  a Deligne cohomology valued characteristic class, we will restrict ourselves
  to differentiable principal $G$-bundles.

  \begin{definition}
  A Deligne characteristic class $d$ (of degree $p$) for principal $G$-bundles with connection is
  an  assignment to any principal $G$-bundle $P$ with connection
  $A$ over $M$ of a  class $d(P, A) \in H^p(M, \cD^p)$
    which is functorial in the sense  that if
  $f \colon N \to M$ then
  $$
  d(f^*(P), f^*(A)) =
  f^*(d(P, A)),
  $$
  where $f^*(P)$ is the pull-back principal $G$-bundle with
  the pull-back connection $f^*(A)$.
  \end{definition}

  Note that if we add two Deligne characteristic classes
  using the group structure in  $H^p(M, \cD^p)$ the
  result is another  Deligne characteristic class.
  Denote by $\cD_p(G)$ the group of all Deligne
  characteristic classes  of degree  $p$ for principal $G$-bundles.

    If $d\in \cD_p(G)$ is a Deligne characteristic class for principal $G$-bundles
     and $P \to M$ is  a principal $G$-bundle, we can
    choose a connection $A$ on $P$, then $d (P, A) \in H^p(M, \cD^p)$.
    Composing with the characteristic class map for Deligne cohomology
    $$c: H^p(M, \cD^p) \to H^{p+1}(M, \ZZ),$$
    we get
    \ba\label{c:d}
    c(d (P, A) )= c\circ d (P, A) \in H^{p+1}(M, \ZZ).
    \na

    \begin{lemma} \label{DC2BG}
    The above map (\ref{c:d}) defines a homomorphism
    $\cD_p(G) \to H^{p+1}(BG, \ZZ)$.
    \end{lemma}
    \begin{proof}
    If we can show that (\ref{c:d})
     is independent of the choice of connections
    then we have defined a characteristic class for $P$,
    which corresponds to  an element in $H^{p+1}(BG, \ZZ)$.
    Here we approximate $BG$ by finite dimensional smooth models (see \cite{NarRam}),
    and use the fact that $H^*(BG, \ZZ)$ is the inductive limit of the cohomology of these 
    models.  
    To see that $c\circ d (P, A)$ doesn't depend of the choice of connections,  
    let $A_0$ and $A_1$ be two
    connections on $P$ and consider the connection $\AA$ on $\hat P = P
  \times \RR \to M \times \RR$
    given by
    \[
    \AA = (1-t)A_0 + t A_1 .
    \]
     Let $\iota_t:  M \to M   \times \RR$ be the
    inclusion map $\iota_t(m) = (m, t)$. It is well known that the induced maps
    on cohomology $\iota_t^* \colon H^p(M \times \RR, \ZZ) \to H^p(M, \ZZ)$
    are all equal and isomorphisms. Moreover $\iota_0^*(\hat P, \AA) = (P, A_0)$
    and $\iota_1^*(\hat P, \AA ) = (P, A_1)$ which imply
    \[\begin{array}{lll}
    c \circ d(P, A_0) &=&
    c \circ d(\iota_0^*(\hat P, \AA )) = \iota_0^*c \circ d(\hat P, \AA ) \\
    &=&  \iota_1^*c \circ  d(\hat P, \AA ) 
    =c \circ d(\iota_1^*(\hat P, \AA ))\\  &  =  &c \circ  d(P, A_1) \in H^{p+1}(M, \ZZ).
  \end{array}
    \]
    Hence we have defined a homomorphism
    $$
    \cD_p(G) \to H^{p+1}(BG, \ZZ).
    $$
    \end{proof}

  Define $I^k(G)$ to be the ring  of invariant polynomials on  the
  Lie algebra of $G$. Then we  have the Chern-Weil homomorphism:
  $$
  cw \colon I^k(G) \to H^{2k}(BG, \RR).
  $$
  If $G$ is compact then this is an isomorphism.
  Define
  $$
  A^{2k}(G, \ZZ) = \{ (\Phi, \phi) \in I^k(G) \times H^{2k}(BG, \ZZ) \mid
  cw(\Phi) = r(\phi) \}.
  $$

  In \cite{CheSim} Cheeger and Simons show that each $(\Phi, \phi) \in
  A^{2k}(G, \ZZ)$ defines a differential character valued
   characteristic  class of degree
  $2k-1$,  whose value on a principal $G$-bundle $P$ over $M$ with connection $A$   is
  denoted by  $$S_{\Phi, \phi}(P, A) \in \check{H}^{2k-1}(M, U(1))$$ (Cf. Remark
  \ref{S:P-A}). Let $c_{\Phi, \phi}(P, A)$ be the element  in $H^{2k-1}(M, \cD^{2k-1})$
  such that under the natural isomorphism
   \[
   \begin{array}{ccc}
   \check{H}^{2k-1}(M, U(1))
  &\to& H^{2k-1}(M, \cD^{2k-1})\\[2mm]
  S_{\Phi, \phi}(P, A)&\mapsto& c_{\Phi, \phi}(P, A).
  \end{array}
  \]
  For the category of principal $G$-bundles with connection whose
  morphisms are connection preserving bundle morphisms, then $c_{\Phi, \phi}(P, A)$
   is a Deligne characteristic class, which is a functorial
  lifting of
  \[
  (\Phi (\disp{\frac{i}{2\pi}}F_A),\phi(P)) \in
  \Omega^{2k}_{cl, 0}(M)\times H^{2k}(M, \ZZ)
  \]
  where $F_A$ is the curvature 2-form of the connection $A$ and $\phi(P)$ is
  the corresponding characteristic class of $P$.
  This defines a map: $ A^{2k}(G, \ZZ)  \to \cD_{2k-1}(G)$.

  In particular, if $G$ is compact and $\phi \in H^{2k}(BG, \ZZ)$ then
  $$\Phi = cw^{-1}(r(\phi)) \in I^k(G)$$ satisfies $cw(\Phi) = r(\phi)$, which means
  $(\Phi, \phi) \in  A^{2k}(G, \ZZ)$. As $\Phi$ is determined by $\phi$
  in this case we write $c_{\phi}\equiv c_{\Phi, \phi}$.  So
  we have a composed map
  $$
  H^{2k}(BG, \ZZ) \to A^{2k}(G, \ZZ)  \to \cD_{2k-1}(G)
  $$
  which sends $\phi$ to $c_\phi$.

    \begin{proposition}\label{DC=BG}
    For a compact Lie group $G$, each element $\phi$ in $ H^{2k}(BG, \ZZ)$ defines
    a degree $2k-1$ Deligne characteristic class $c_\phi$
      in $\cD_{2k-1}(G)$  such that $c_\phi\mapsto \phi$ under the homomorphism in Lemma
      \ref{DC2BG}.
    \end{proposition}
        \begin{proof} From the above discussion and
    the definition of the Deligne characteristic class, we obtain that,
    given a principal $G$-bundle $P$ with a connection $A$ over $M$, 
    $c_\phi(P, A)= c_{\Phi, \phi}(P, A)$ is a functorial  lifting of
  \[
  (\Phi (\disp{\frac{i}{2\pi}}F_A),\phi(P)) \in
  \Omega^{2k}_{cl, 0}(M)\times H^{2k}(M, \ZZ).
  \]
  This implies that the corresponding characteristic class of $P$ is $\phi(P)$. From
  the  bijective correspondence between degree $2k$ 
  characteristic classes of principal
  $G$-bundles  and elements of $H^{2k}(BG, \ZZ)$, we know
  that $c_\phi\mapsto \phi$ under the homomorphism in Lemma
      \ref{DC2BG}.
      \end{proof}

   Fixing a smooth infinite dimensional model of
    $EG\to BG$ by embedding the compact semi-simple Lie group
    $G$ into $U(N)$ and letting $EG$ be the Stiefel manifold of
    $N$ orthonormal vectors in a separable complex Hilbert space, we know
    that  the Deligne cohomology group $H^{2k-1}(BG,\cD^{2k-1})$
    is well-defined.
    
   Let $\AA$ be    a universal connection on $EG$, $\phi \in H^{2k}(BG, \ZZ)$
    defines a degree $2k-1$ Deligne characteristic class
    \[
    c_ {\Phi, \phi} (EG, \AA)  \in H^{2k-1}(BG,\cD^{2k-1}),
    \]
    where $\Phi \in I^{k}(G)$ satisfies $cw(\Phi) = r(\phi)$. Then the
    following commutative diagram
    \[
 \xymatrix{
 H^{2k-1}(BG, \cD^{2k-1}) \ar[d] \ar[r]& H^{2k}(BG, \ZZ)\ar[d]\\
 \Omega^{2k}_{cl, 0}(BG) \ar[r] & H^{2k}(BG, \RR)}
 \]
     shows that the map $\phi \mapsto c_ \phi (EG, \AA)$ refines
   the Chern-Weil homomorphism.

  \section{From Chern-Simons to Wess-Zumino-Witten}\label{section:CS2WZW}

Let $G$ be  a compact, connected, semi-simple Lie group. 
  In \cite{DijWit}  Dijkgraaf and Witten discuss a correspondence map
  between three dimensional Chern-Simons gauge
   theories and Wess-Zumino-Witten models associated to the compact Lie group
   $G$  from the topological actions viewpoint,
which naturally involves the transgression map
  \[
  \tau:  H^4(BG, \ZZ) \to H^3(G, \ZZ).
  \]
  (To be precise, $\tau$ is actually
the inverse of the transgression in Borel's study of
  topology of Lie groups and characteristic classes but this is of no
real importance.)

  We recall the definition of $\tau$.
 We take a class $\phi \in H^4(BG, \ZZ)$
  and pull its representative $\phi$  back to $\pi^*(\phi) $, a four-cocycle 
  on $EG$. As $EG$ is contractible we have that $\phi = d\tau_\phi$
  for a three-cocycle $\tau_\phi$ on $EG$. Restricting $\tau_\phi$ to a fibre which
  we identify with $G$ it is easy to show that the result
  is a closed cocycle defining an element of $H^3(G, \ZZ)$ and
  that moreover this cohomology class is independent of all 
  choices made.  

 It is shown in \cite{DijWit} that three dimensional Chern-Simons gauge
   theories with gauge group $G$ can be classified by the integer
cohomology group    $H^4(BG, \ZZ)$, and
  conformally invariant sigma models in two dimension with 
  target space a compact Lie group
  (Wess-Zumino-Witten models)  can be classified by $H^3(G, \ZZ)$.
  Recall the   commutative diagram
    \[
 \xymatrix{
 H^{3}(BG, \cD^{3}) \ar[d] \ar[r]& H^{4}(BG, \ZZ)\ar[d]\\
 \Omega^{4}_{cl, 0}(BG) \ar[r] & H^{4}(BG, \RR).}
 \]
  To classify the exponentiated Chern-Simons action  in  three dimensional
   Chern-Simons gauge theories, we propose the following mathematical
   definitions of a three dimensional  Chern-Simons gauge theory and
   a Wess-Zumino-Witten model.
  
  \begin{definition} We make the following definitions: 
  \begin{enumerate}
  \item  A three dimensional {\em Chern-Simons gauge theory}
with gauge group
  $G$ is defined to be a Deligne   characteristic class of degree $3$ for
  a principal $G$-bundle with connection.
We denote the group of all three dimensional Chern-Simons gauge theories
   with gauge group   $G$  by   $CS(G)$.
 \item  A  {\sl Wess-Zumino-Witten model} on $G$ is defined to be
   a Deligne class on $G$ of degree $2$ and we
  denote the group of all such  by $WZW(G)$.
  \end{enumerate}
  In brief, $CS(G) = \cD_3(G)$ and $WZW(G) = H^2(G, \cD^2)$.
 \end{definition}

  With these preliminaries taken care of we can now explain our
  refined  geometric definition of
  the Dijkgraaf-Witten map and discuss its image.
  Firstly, we give a more geometric definition of the
  transgression map $\tau$,  which has the advantage that it can
  be lifted to define a correspondence  map  from $CS(G)$ to $WZW(G)$.
   We do this by constructing a canonical
  $G$-bundle $\cP$ with connection $\mathbb A$ on the manifold
  $S^1 \times G$. It follows that if $d \in CS(G)$ then
  $d(\cP, {\mathbb A}) \in H^3(S^1 \times G , \cD^3)$. We can
  integrate along $S^1$  with the Deligne characteristic class
  in  $ H^3(S^1 \times G , \cD^3)$ to   get the required map
  $$
  \int_{S^1} d(\cP, \mathbb A) \in H^2(G, \cD^2) = WZW(G).
  $$

  We want to show that for any $G$ there is a natural
  $G$ bundle over $S^1 \times G$ with connection. To do this
  it is convenient to work with {\em pre-$G$-bundles}.

 \begin{definition} 
  A {\sl pre-$G$-bundle} is a pair $(Y, \hat{g})$ where $\pi: Y \to M$
  is a surjective submersion and    $\hat{g} \colon Y^{[2]} \to G$ such that
  $$
  \hat{g}(y_1, y_3)  = \hat{g}(y_1, y_2)\hat{g}(y_2, y_3)
  $$
  for any $y_1, y_2, y_3$ all in the same fibre of $\pi
  \colon Y \to M$.  Here we denote by $Y^{[p]}$ the $p$-fold
  fibre product of $\pi\colon Y \to M$.
  \end{definition}

  If $P \to M$ is a principal $G$ bundle, there is a canonical
  map $\hat{g} \colon P^{[2]} \to G$ defined by 
  \[
  p_1 \hat{g}(p_1, p_2) = p_2.
  \]
  Then $(P, \hat{g})$ is a pre-$G$-bundle. Conversely
  if $(Y, \hat{g})$ is a pre-$G$-bundle over $M$,
   we can construct a principal $G$-bundle
  $P \to M$ as follows. Take $Y \times G$ and define 
  \[
  (y_1, h_1)\sim  (y_2, h_2)\quad  \text{if} \  \pi(y_1) = \pi(y_2)\ 
  \mathrm{and}\ h_1 \hat{g}(y_1, y_2) = h_2.
  \] 
  The space of equivalence classes $P= Y \times G / \sim $ is a principal
 $G$-bundle over $M$  with right $G$-action on equivalence classes $[y, g]$
given by $[y, h]\cdot g = [y, hg]$

  Two pre-$G$-bundles $(Y, \hat{g}_1)$ and $(X, \hat{g}_2)$ give rise to isomorphic
  principal $G$ bundles if and only if there is an $\hat{h} \colon Y
  \times_\pi X \to G$
  such that
  $$
  \hat{h}(y_1, x_1) \hat{g}_1(y_1, y_2) = \hat{g}_2(x_1, x_2) \hat{h}(y_2, x_2)
  $$
  for any collection of points $y_1, y_2 \in Y$, $x_1, x_2 \in X$ mapping
  to the same point in $M$. A pre-$G$-bundle $(Y, \hat{g})$ is {\sl trivial}
  if there is an $\hat{h} \colon Y \to G$ such that 
  \[
  \hat{g}(y_1, y_2) = \hat{h}(y_1)^{-1}\hat{h}(y_2)
  \]
  for every $(y_1, y_2) \in Y^{[2]}$.

 Given  a pre-$G$-bundle $(Y, \hat{g})$ over $M$,  we denote by
  $\hat{g}^{-1}d\hat{g}$ the
  pull-back by $\hat{g} \colon Y^{[2]} \to G$ of the Maurer-Cartan form.
  Then $\hat{g}^{-1}d\hat{g}$ is a $\g$ (the Lie algebra of $G$)-valued . 
  Let $A$ be a $\g$-valued   one-form on $Y$. We say that $A$ is a connection
  for the pre-$G$-bundle $(Y, \hat{g})$  if
  $$
  \pi_1^*(A) = ad (\hat{g}) \pi_2^*(A)  - \hat{g}^{-1} d\hat{g}
  $$
  where  $\pi_1, \pi_2 \colon Y^{[2]} \to Y$ are
  the projections and we denote the adjoint action of $ G$
  on its Lie algebra by $ad(\hat{g})$.  It is easy to check that there
  is a one-to-one correspondence between connections on a
  pre-$G$-bundle and connections on the associated principal $G$ bundle.

  We wish to define a $G$ bundle on $S^1 \times G$ with
  connection. From the previous discussion it suffices to define a
  pre-$G$-bundle with connection.
  Let  $\cA$  be all smooth maps $h$ from $\RR$ to $G$
  with $h^{-1}dh$
  periodic and $h(0) = 1$.  
  Define $\pi \colon \cA \to G$ by
  $\pi(h) = h(1)$.  Notice that if $\pi(g)
  = \pi(h)$ then $g = h\c$   where $\c$ is a smooth map from $[0,1]$ to $G$
  with $\c^{-1}d\c$
  periodic and $\c(1) =1 = \c(0)$.  Such a $\c$ is actually
  a smooth
  based loop in the based loop group $\Omega G$. We can identify
  $\cA$ with the space of $G$-connections on the circle $S^1= \RR/\ZZ$.
  Then  $\cA$   is contractible and $\pi \colon \cA \to G$ is the holonomy map.  
   Hence,   $\pi \colon \cA \to G$ is a universal $\Omega G$-bundle, and $G$
  is a classifying space $B\Omega G$ of $\Omega G$.
    Let $Y = \cA \times S^1 \to G \times S^1$.
  Define $\hat{g}  \colon Y^{[2]} \to G$ by 
  \[
  \hat{g}(h_1, h_2, \theta) = h_1(\theta)^{-1}
  h_2(\theta).
  \]
   Then the pair $(Y,\hat{g})$ is a pre-$G$-bundle over $G \times S^1$.
  Let $\hat Y = \RR \times \cA$ and the projection $\hat Y \to Y $
  induced by $\RR \to S^1$.  Define $\hat h \colon \hat Y  \to G$ by
  $\hat{h}(t, h) = h(t)$.
   $\hat{h}^{-1}d\hat{h}$ being periodic descends to a $\g$-valued 
  one-form on $Y$. It is straightforward to check that
  this defines  a connection $A$ for the pre-$G$-bundle $(Y, \hat{g})$ over $G\times S^1$.

The principal $G$-bundle  over $G \times S^1$ corresponding to the 
pre-$G$-bundle $Y= \cA\times S^1 \to G\times S^1$
can be obtained as follows (Cf. \cite{CarMic} and \cite{MurSte}). Denote by
\[
\cP = \disp{\frac{\cA \times S^1 \times G}{\Omega G}},
\]
the quotient space of $\Omega G$-action on $\cA\times S^1 \times G$, where the
$\Omega G$-action is given by, for $\gamma \in \Omega G$ and
$(h, \theta, g)\in \cA \times S^1 \times G$,
\ba\label{action}
\gamma \cdot (h, \theta, g)= (h\gamma, \theta, \gamma (\theta)^{-1} g).
\na
Notice that $\cP$ admits a natural free $G$-action from the right multiplication
on $G$-factor.  The connection $A$ on the pre-$G$-bundle $(Y,\hat{g})$ defines a
natural connection $\mathbb A$ on $\cP$.

\begin{definition}\label{CS2WZW:map}
The canonical principal $G$-bundle over $G\times S^1$ is given by $\cP$ with
connection $\mathbb A$. The correspondence map from three dimensional
Chern-Simons gauge
theories $CS(G)$ to Wess-Zumino-Witten models $WZW(G)$ is defined to be
\[
CS(G) = \cD_3(G) \ni d \mapsto \disp{\int_{S^1}} d(\cP, \mathbb A) \in H^2(G, \cD^2)
=WZW(G).
\]
Denote this map by
$\Psi: CS(G) \longrightarrow WZW(G).$
\end{definition}

The next proposition shows that the map $\Psi$
descends to the natural transgression map from
$H^4(BG, \ZZ)$ to $H^3(G, \ZZ)$ and
 hence $\Psi$ refines the
Dijkgraaf-Witten correspondence.

\begin{proposition}
\label{true:trans}
The correspondence map from $CS(G)$ to $WZW(G)$
induces the natural transgression map $\tau$: $H^4(BG, \ZZ)\to  H^3(G, \ZZ)$.
\end{proposition}
\begin{proof}
We first give another construction of $\tau$. Let
$EG\to BG$ be the universal $G$-bundle, then (as is well 
known) the $\Omega G$ bundle 
\[
\tilde{\pi}: \Omega EG \to \Omega BG
\]
formed by applying the based loop functor to $EG\to BG$ 
gives another model of the universal $\Omega G$-bundle.  
In particular we have a homotopy equivalence $\Omega BG 
\stackrel{\simeq}{\to} G$ which lifts to an $\Omega G$-equivariant homotopy 
equivalence $\Omega EG\stackrel{\simeq}{\to} \cA$.  
This leads to the
isomorphism:
\[
H^3(\Omega BG, \ZZ) \to H^3(G, \ZZ).
\]
On the other hand, the natural evaluation map:
\[
ev: \Omega BG \times S^1 \longrightarrow BG
\]
defines a pull-back map
\[
ev^*:  H^4(BG, \ZZ) \to H^4( \Omega BG \times S^1, \ZZ),
\]
from which the integration along $S^1$ gives rise to another
construction of the transgression map:
\ba\label{trans}
\disp{\int}\circ ev^*:  H^4(BG, \ZZ) \to H^3(\Omega BG, \ZZ)  \cong H^3(G, \ZZ).
\na

{}From the homotopy equivalence between the two universal $\Omega G$-bundles:
$\pi: \cA \to G$ and $\tilde{\pi}: \Omega EG \to \Omega BG$, we get
a homotopy equivalence of two principal $G$-bundles:
\ba\label{homotopy}
\cP = \disp{\frac{\cA \times S^1 \times G}{\Omega G}} \sim
\disp{\frac{ \Omega EG\times S^1 \times G}{\Omega G}}.
\na
Here the $\Omega G$ action on $\Omega EG\times S^1 \times G$ is given by
the similar action on $ \cA \times S^1 \times G$ as in (\ref{action}).
$\frac{ \Omega EG\times S^1 \times G}{\Omega G}$ is a principal $G$-bundle
over $\Omega BG\times S^1$.

Now we show that the pull-back of the universal $G$-bundle: $EG\to BG$ via
the evaluation map $ev$, which is
\ba\label{ev}
ev^*(EG) = \bigl( \Omega BG \times S^1\bigr)\times_{BG} EG
\na
 is isomorphic to
$\frac{ \Omega EG\times S^1 \times G}{\Omega G}$ as principal $G$-bundles.
The isomorphism map
\ba\label{=}
\disp{\frac{ \Omega EG\times S^1 \times G}{\Omega G}}
\to \bigl( \Omega BG \times S^1\bigr)\times_{BG} EG
\na
is given by
\[
[(\tilde \c, \theta, g)] \mapsto [ (\tilde{\pi}(\tilde{\c}), \theta,
\tilde{\c}(\theta)\cdot g)].
\]
Here $(\tilde \c, \theta, g) \in \Omega EG\times S^1 \times G$, $\tilde\pi$ is the
map $\Omega EG \to \Omega BG$,  $\tilde{\c}(\theta)$ is the image of
the evaluation map on $\Omega EG\times S^1$ and the action of $g$
on $\tilde{\c}(\theta)$ is induced from the right $G$-action on
the universal $G$: $EG\to BG$.

It is easy to check that
(\ref{=}) is a well-defined $G$-bundle isomorphism by direct calculation:
\[\begin{array}{lll}
[(\tilde {\c}\cdot \c , \theta, \c(\theta)^{-1}g)]  &\mapsto &
 [(\tilde{\pi}(\tilde{\c}\cdot \c), \theta,
(\tilde{\c}\gamma)(\theta)\cdot(\c(\theta)^{-1} g))]\\[2mm]
&=& [ (\tilde{\pi}(\tilde{\c}), \theta,
\tilde{\c}(\theta)\cdot g)],
\end{array}
\]
and for $g'\in G$
\[
[(\tilde \c, \theta, g g')]\mapsto [ (\tilde{\pi}(\tilde{\c}), \theta,
\tilde{\c}(\theta)\cdot gg')].
\]

Then from  (\ref{homotopy}), (\ref{ev}) and (\ref{=}), we obtain
a homotopy equivalence of two principal $G$-bundles:
\ba\label{homotopy:1}
\begin{array}{ccc}
\cP & &
ev^*(EG)\\[2mm]
\downarrow &\sim& \downarrow\\[2mm]
S^1\times G && S^1\times \Omega BG.
\end{array}
\na
Hence, with the definition of the integration map (\ref{integration:1})
given by (\ref {integration:2}) and (\ref {integration:3}), we see that
 the transgression map $\tau= \int\circ ev^*:
 H^4(BG, \ZZ) \to H^3(G, \ZZ)$  agrees with the
map induced by our correspondence $\Psi:CS(G)\to WZW(G)$.
  \end{proof}

 For a compact Lie group, Proposition \ref{DC=BG} tells us that there exists a one-to-one
 map
 \[
 CS(G) \to H^3(BG, \cD^3)\]
 and the exact sequence (\ref{exact:1}) implies the exact sequence
 \[
  0 \rightarrow
  \Omega^2(G)/ \Omega^2_{cl, 0}(G) \rightarrow WZW(G)) \rightarrow H^{3}(G, \ZZ)
  \rightarrow 0.
  \]
  As in general, the map $\tau: H^4(BG, \ZZ) \to H^{3}(G, \ZZ)$ is not
  surjective, the map
  \[
 \Psi:   CS(G) \to WZW(G)
  \]
  is not surjective either, for a general compact semi-simple Lie group.
We will see that the
  Wess-Zumino-Witten models from  the image of $\Psi$
  exhibit some special properties by exploiting bundle gerbe theory.

We give a summary of how various bundle gerbes enter. First it is now
well understood how,
 given a WZW model, we can define an associated  bundle gerbe
  over the group $G$, as $WZW(G) = H^2(G, \cD^2)$ is the space of stable
  isomorphism classes of bundle gerbes with connection and curving 
  \cite{CMM,MurSte0}.
We will see in the   Section \ref{b2g} that $H^4(BG, \ZZ)$  is the space
  of stable equivalence classes of  bundle $2$-gerbes on $BG$.
Thus an element in
   $H^4(BG, \ZZ)$ defines a class of bundle $2$-gerbes on $M$ associated to
  a principal $G$-bundle over $M$ using the pullback construction of the
  classifying map.
The corresponding transgressed element in
  $H^3(G, \ZZ)$ defines  a bundle gerbe over $G$. In fact, as
  $H^3(BG, \RR) =0$, we know that  the third Deligne cohomology group
  $H^3(BG, \cD^3)$ is determined by the  following commutative diagram:
  \[
  \begin{array}{ccc}
 H^3(BG, \cD^3) &\sta{c}{\longrightarrow}& H^4(BG, \ZZ)\\[2mm]
  \downarrow^{curv}&&\downarrow\\[2mm]
  \Omega^4_{cl, 0}(BG)&\longrightarrow & H^4(BG, \RR).
  \end{array}
  \]
  In the last Section we showed that
given an element $ \phi \in H^4(BG, \ZZ)$,  we can take a
  $G$-invariant polynomial $\Phi \in I^2(G)$ corresponding to $\phi$, then
  for a connection $\AA$ on the universal bundle $EG$ over $BG$,
  \[
  (\phi, \Phi(\disp{\frac{i}{2\pi }}F_\AA)) \in H^4(BG, \ZZ)\times_{ H^4(BG, \RR)}
  \Omega^4_{cl, 0}(BG)
  \]
  represents a degree $3$ Deligne class
  \[
  c_\phi(EG, \AA) \in  H^3(BG, \cD^3).
  \]
In a following Section we will define a
universal Chern-Simons bundle 2-gerbe determined by $c_\phi(EG, \AA)$.
We will then pull it back to define,
for a principal $G$-bundle $P$ with
connection $A$ over $M$, a
  Chern-Simons bundle 2-gerbe over $M$, with 2-curving given by
  the Chern-Simons  form  associated to $(P, A)$.

  \begin{remark}
 Brylinski   \cite{Bry} defines a generalisation of $H^3(G,\ZZ)$, called the
  differentiable cohomology, denoted $H^3_{\mathit{diff}}(G,U(1))$ for which there is
  an isomorphism $H^3_{\mathit{diff}}(G,U(1)) \cong H^4(BG,\ZZ)$. In
  low dimensions these differentiable cohomology classes have the following
  interpretations \cite{Bry}:
  \begin{center}
  \begin{tabular}{ccl}
  $H^1_{\mathit{diff}}(G,U(1))$ & $\cong$ &  smooth homomorphisms $G \rightarrow U(1)$\\
  $H^2_{\mathit{diff}}(G,U(1))$ & $\cong$ & isomorphism classes of central extensions of
  $G$ by $U(1)$\\
  $H^3_{\mathit{diff}}(G,U(1))$ & $\cong$ & equivalence classes of multiplicative
  $U(1)$-gerbes on $G$.
  \end{tabular}
  \end{center}
  His multiplicative $U(1)$-gerbes motivated
  us to define multiplicative bundle gerbes.
  \end{remark}

  \section{Bundle 2-gerbes} \label{b2g}

 Bundle 2-gerbe theory on $M$  is developed in \cite {Ste}.
A bundle 2-gerbe
 with connection and curving, defines a degree 3 Deligne class
in $H^3(M, \cD^3)$.  In \cite{Joh}
 it is shown that the group of stable equivalence classes of bundle 2-gerbes
 with connection and curving is isomorphic to $H^3(M, \cD^3)$.
We review the results in
 \cite{Ste} and \cite{Joh},
then we define multiplicative bundle gerbes
 on  compact Lie group.

  We begin with the definition of a simplicial bundle gerbe
  as in  \cite{Ste} on   a simplicial manifold $X_\bullet =\{X_n\}_{n\ge 0}$ with face operators
  $d_i: X_{n+1} \to X_{n}$ ($i= 0, 1, \cdots, n+1$).  
  We remark that the simplicial manifolds
  we use in this paper are not required to have degeneracy operators (see
  \cite{Dup}).

  \begin{definition}  \label{simplicial:gerbe} (Cf. \cite{Ste})
  A simplicial bundle gerbe on a simplicial manifold $X_\bullet$ consists of the following data.
  \begin{enumerate}
  \item A bundle gerbe $\cG$ over $X_1$.
  \item A bundle gerbe stable isomorphism $m: d_0^*\cG \otimes d_2^* \cG  \to d_1^*\cG$
  over $X_2$, where $d_i^*\cG$ is the pull-back bundle gerbe over $X_2$.
  \item The bundle gerbe stable isomorphism $m$ is associative up to a natural transformation,
  called an associator,
  \[
  \phi: \qquad   d_2^*m \circ ( d_0^*m \otimes Id)
    \to   d_1^*m \circ (Id \otimes d_3^*m)
  \]
  between the induced stable isomorphisms  of bundle gerbes over $X_3$.
  The line bundle  $\cL_\phi$ over $X_3$ induced  by $\phi$ 
  admits a trivialisation section $s$ such that
    $\delta (s)$ agrees with the canonical trivialisation of 
   \[ 
   d_0^*\cL_\phi \otimes d_1^*\cL_\phi^* \otimes d_2^*\cL_\phi
   \otimes d_3^*\cL_\phi^*  \otimes d_4^*\cL_\phi.
   \]
  \end{enumerate}
  If in addition the bundle gerbe $\cG$ is equipped with a connection and a curving,
  and $m$ is a stable isomorphism of bundle gerbes with connection and curving,
  we call it a simplicial bundle gerbe with connection and curving on $X_\bullet$.
  \end{definition}

  \begin{remark}
  \label{some:mercy} For those unfamiliar with 2-gerbes we offer
the following amplification.
  \begin{enumerate} \item A bundle gerbe stable isomorphism $m$
  in Definition \ref{simplicial:gerbe}
  is a fixed trivialisation of the bundle gerbe
  $$\delta(\cG) =d_0^*\cG \otimes d_1^*\cG^* \otimes d_2^* \cG$$
  over $X_2$, 
  where $d_1^*\cG^*$ is the dual bundle gerbe of $d_1^*\cG$.
  (See \cite{Mur}\cite{MurSte0}
   for  various operations on bundle gerbes and the
   definition of a bundle gerbe stable isomorphism.)
   \item With the understanding of  $m$ in Definition \ref{simplicial:gerbe}
    as a fixed  trivialisation of the bundle gerbe 
    $d_0^*\cG \otimes d_1^*\cG^* \otimes d_2^* \cG$
  over $X_2$ , we can see
    that     $d_2^*m \circ ( d_0^*m \otimes Id)$ and $d_1^*m \circ (Id \otimes d_3^*m) $
    represent two trivialisations of the bundle gerbe   over $X_3$.
 This induces a line bundle over $X_3$ (Cf. \cite{Mur}), called
the $\phi$-induced or associator line bundle. A simplicial bundle gerbe
$\cG$ requires that this associator line bundle is trivial and the
trivialisation section $s$ satisfies the natural cocycle condition. 
  \item A simplicial bundle gerbe with connection and curving has
in its definition, a restrictive
  condition, as it requires that the bundle gerbe stable isomorphism $m$
  preserves connections and curvings. This implies,
  \ba\label{curv:0}
  d_0^*(curv(\cG)) - d_1^*(curv(\cG)) + d_2^*(curv(\cG)) = 0.
  \na
   For the simplicial bundle gerbe constructed in this
   paper, the underlying bundle gerbe is often equipped with a connection and
  curving, but we shall not require that
   the bundle gerbe stable isomorphism $m$ preserves the
 connection and curving.
  In section \ref{multi:WZW}, instead of (\ref{curv:0}), we have
   \[
    d_0^*(curv(\cG)) - d_1^*(curv(\cG)) + d_2^*(curv(\cG)) = dB,
    \]
    for some 2-form $B$ on $X_2$, which is not necessarily obtained from
    a 2-form on $X_1$ by the $\delta$-map $d_0^*- d_1^*+ d_2^*$.
    \end{enumerate}
  \end{remark}

  For a smooth submersion $\pi: X\to M$, there is a natural associated simplicial
 manifold $X_\bullet = \{X_n\}$ (which one might well think of as the `nerve'
associated to $\pi:X\to M$) with $X_n$ given by
  \[
  X_n = X^{[n+1]} = X\times_M X\times_M \cdots \times_M X
  \]
  the $(n+1)$-fold fiber product of $\pi$, and face operators
  $d_i = \pi_{i+1}: X_{n} \to X_{n-1}$ ($i= 0, 1, \cdots, n$) 
  are given by the natural projections
  from $X^{[n+1]}$ to $X^{[n]}$ by omitting the entry in $i$ position for $\pi_i$.
  For an exception,  we denote by 
  $EG^{[\bullet]}$ the  associated simplicial
 manifold $\{EG^{[n]}\}$ for the universal bundle $\pi: EG \to BG$,.

  Now we recall the definition of bundle 2-gerbe on a smooth manifold $M$ from
  \cite{Ste} and \cite{Joh}.
  
  \begin{definition}
  A bundle 2-gerbe on $M$ consists of a quadruple of smooth manifolds
  $(\cQ, Y; X, M)$ where $\pi: X\to M$ is a smooth, surjective submersion, and
  $(\cQ, Y; X^{[2]})$ is a simplicial bundle gerbe on the simplicial manifold
  $X_\bullet = \{X_n = X^{[n+1]}\}$ associated to $\pi: X\to M$.
  \end{definition}

It is sometimes convenient to describe bundle $2$-gerbes 
using the language of $2$-\emph{categories} (see for example 
\cite{KellyStreet}).  One first observes that transformations 
between stable isomorphisms provide $2$-morphisms making 
the category $\mathbf{BGrb}_M$ of bundle gerbes over $M$
and stable isomorphisms between bundle gerbes into a weak $2$-category or 
\emph{bi-category} (Cf. \cite{Ste}). Note that the space of 2-morphisms between
two stable isomorphisms is one-to-one corresponding to
the space of line bundles over $M$.

Consider the face operators $\pi_i: X^{[n]} \to X^{[n-1]}$ on the simplicial 
manifold $X_\bullet = \{X_n = X^{[n+1]}\}$. We can define a bifunctor
\[
\pi_i^*:   \mathbf{BGrb}_{X^{[n-1]}}  \longrightarrow
\mathbf{BGrb}_{X^{[n]}} 
\]
sending objects, stable isomorphisms and 2-morphisms to 
the pull-backs by $\pi_i$ ($i=1,\dots, n$).
 One can then use this language to 
describe the data of a bundle $2$-gerbe as follows.  
A bundle $2$-gerbe on $M$ consists of the data of a 
smooth surjective submersion $\pi\colon X\to M$ together 
with
\begin{enumerate}
\item  An object $(\cQ,Y,X^{[2]})$
in $\mathbf{BGrb}_{X^{[2]}}$.
\item  A stable isomorphism $m$: $\pi_1^*\cQ\otimes \pi_3^*\cQ  
\to \pi_2^*\cQ$ in $\mathbf{BGrb}_{X^{[3]}}$ defining the bundle 2-gerbe
product which is associative up to a 2-morphism $\phi$
in $\mathbf{BGrb}_{X^{[4]}}$. 
\item The 2-morphism $\phi$ satisfies a natural coherency
condition in $\mathbf{BGrb}_{X^{[5]}}$.
\end{enumerate}

We now briefly pause to describe some new notation 
which provides a good way to encode the simplicial bundle gerbe 
data (Cf. Definition \ref{simplicial:gerbe} and Remark \ref{some:mercy}).
  We  define 
maps $\pi_{ij}\colon X^{[n]}\to X^{[2]}$ for $n>2$ which 
send a point $(x_1,\ldots,x_{n})$ of $X^{[n]}$ 
to the point $(x_i,x_j)\in X^{[2]}$.   
It is clear that these maps 
can be written (non-uniquely) in terms of the $\pi_i$ 
(the non-uniqueness stems from the simplicial identities 
satisfied by the face maps $\pi_i$'s).  Let us write $\cQ_{ij}$ 
for $\pi_{ij}^*\cQ$.  For example, the bundle gerbe $\cQ_{12}$ over $X^{[3]}$
 is the pull-back $\pi_3^*\cQ$ of $\cQ$.

Returning to the definition of bundle $2$-gerbe, the next 
part of the definition requires that there is a stable isomorphism 
$m: \cQ_{23}\otimes \cQ_{12} \to \cQ_{13}$ of bundle gerbes 
over $X^{[3]}$ together with a natural transformation called an associator
  $$
   \phi:\  \pi_3^*m \circ ( \pi_1^*m \otimes Id)
    \to \pi_2^*m \circ (Id \otimes  \pi_4^*m)
  $$
   which is a 2-morphism in the bi-category $\mathbf{BGrb}_{X^{[4]}}$,
  between the induced stable isomorphisms of bundle gerbes over $X^{[4]}$
 making the following diagram commute up to an associator $\phi$
 (represented by a 2-arrow in the diagram): 
\ba
\xymatrix{ \label{2-arrow} 
&\cQ_{34}\otimes \cQ_{23}\otimes \cQ_{12} \ar[d]_{ \pi_1^*m \otimes Id}
\ar[r]^-{Id \otimes \pi_4^*m} & \cQ_{34}\otimes \cQ_{13} \ar@2{->}[dl]_-{\phi} 
\ar[d]^{\pi_2^*m} \\ 
&\cQ_{24}\otimes \cQ_{12} \ar[r]^-{\pi_3^*m} & \cQ_{14} } 
\na 
Here we write $\pi_1^*m$ as a stable isomorphism 
$\cQ_{34}\otimes \cQ_{23} \to \cQ_{24}$ over $X^{[4]}$, similarly for
$\pi_2^*m,  \pi_3^*m$ and $ \pi_4^*m$. Hence, the  associator $\phi$
as a 2-morphism in $\mathbf{BGrb}_{X^{[4]}}$ defines a line bundle $\cL_\phi$ over 
$X^{[4]}$, which is required to have a trivialisation section.

In order to write the efficiently 
 coherence condition satisfied by the natural transformation $\phi$, we 
need one last piece of new notation.  
Let us write $\cQ_{ijk} = \cQ_{jk}\otimes \cQ_{ij}$, 
$\cQ_{ijkl} = \cQ_{kl}\otimes \cQ_{jk}\otimes \cQ_{ij}$ 
and so on.  So for example, $\cQ_{123} = \cQ_{23}\otimes \cQ_{12}$,
$ \cQ_{1234}= \cQ_{34}\otimes \cQ_{23}\otimes \cQ_{12}$ and
the diagram (\ref{2-arrow}) in $\mathbf{BGrb}_{X^{[4]}}$ can be written  as
\[
\xymatrix{ 
\cQ_{1234} \ar[d] 
\ar[r] & \cQ_{134} \ar@2{->}[dl]_-{\phi} \ar[d] \\ 
\cQ_{124} \ar[r] & \cQ_{14} } 
\] 
which is commutative if and only if $\phi$ is the
identity 2-morphism denoted by $Id$.

The coherency condition satisfied by the natural 
transformation  $\phi$ can then be viewed from the
 following  two equivalent diagrams in $\mathbf{BGrb}_{X^{[5]}}$
 calculating the associator (2-morphism) between the two induced
 stable isomorphisms from $\cQ_{12345}$ to $\cQ_{15}$ (one is
 $\cQ_{12345}\to \cQ_{1345}\to \cQ_{145}\to \cQ_{15}$, and the other is
 $\cQ_{12345}\to \cQ_{1235}\to \cQ_{125}\to \cQ_{15}$.)
\ba\label{2-arrow:commute}
\xymatrix{       
& \cQ_{12345} \ar[dl] \ar[dr] \ar[rr] & & 
\cQ_{1345} \ar@2{->}[dl]_{\pi_5^*\phi} \ar[dr] & \\ 
\cQ_{1235} \ar[dr] & & \cQ_{1245} \ar[dl] 
\ar[rr] \ar@2{->}[ll]_{\pi_1^*\phi} & & \cQ_{145} \ar[dl] 
\ar@2{->}[dlll]_{\pi_3^*\phi}                            \\ 
& \cQ_{125} \ar[rr] &  & \cQ_{15} & \\ 
& \cQ_{12345} \ar[dl] \ar[rr] & & \cQ_{1345} 
\ar[dl] \ar[dr] \ar@2{->}[dlll]_{Id} &          \\ 
\cQ_{1235} \ar[dr] \ar[rr] & & \cQ_{135} \ar[dr] \ar@2{->}[dl]_{\pi_4^*\phi}  
& & \cQ_{145} \ar@2{->}[ll]_{\pi^*_2\phi} \ar[dl]     \\ 
& \cQ_{125} \ar[rr] & & \cQ_{15}                } 
\na 
which implies the canonical isomorphism  of two
trivial line bundles over $X^{[5]}$
\[
\pi^*_1\cL_\phi \otimes \pi^*_3\cL_\phi \otimes \pi^*_5\cL_\phi
\cong \pi^*_2\cL_\phi \otimes \pi^*_4\cL_\phi.
\]

\begin{remark}  The first example of a bundle 2-gerbe
is the tautological bundle 2-gerbe constructed
  in \cite{CMW} over a 3-connected manifold $M$
  with a closed 4-form $\Theta\in \Omega^4_{cl, 0}(M)$,
see \cite{Ste} for a detailed proof and more examples.
\end{remark}

  \begin{definition}
   Let $(\cQ, Y; X, M)$ be bundle 2-gerbe on $M$. A {\sl bundle 2-gerbe connection}
    on $\cQ$  is a pair $(\nabla, B)$ where $\nabla$ is a bundle gerbe
   connection on the bundle gerbe $(\cQ, Y; X^{[2]})$ and $B$ is a curving
   for the bundle gerbe with connection $(\cQ, \nabla)$,
    whose bundle gerbe curvature
   $\omega$ on $X^{[2]}$ satisfies $\delta (\omega) =0$, where
   $\delta = \pi^*_1- \pi^*_2 + \pi^*_3: \Omega^*(X^{[2]}) \to \Omega^*(X^{[3]})$.
   Then we can solve the equation
   \[
   \omega = (\pi^*_1-\pi^*_2)(C)
   \]
   for a three form $C$ on $X$, such a choice of $C$ is called a {\sl 2-curving }
   for the    bundle 2-gerbe $(\cQ, Y; X, M)$, or simply the {\sl bundle 2-gerbe curving}.
    Then $dC= \pi^*(\Theta)$ for a closed four form
    $\Theta$ on $M$, which is called the {\sl bundle 2-gerbe curvature}.
 \end{definition}

     Locally, as in \cite{Ste} and \cite{Joh},
    a bundle 2-gerbe on $M$ with connection and curving is determined by
    a degree 3 Deligne class
    \ba\label{local:b2g:d-class}
    [(g^{\ }_{ijkl}, A_{ijk}, B_{ij}, C_{i})]\in H^3(M, \cD^3)
    \na
    for a good cover $\{U_i\}$ of $M$, over which there are
local sections $s_i: U_i \to X$.
    Then $C_i =s_i^* C$. Using $(s_i, s_j)$, we can pull-back
    the bundle gerbe $(\cQ, Y, X^{[2]})$ to $U_{ij}= U_i\cap U_j$, such that
    $\cQ_{ij} = (s_i, s_j)^*\cQ$ is trivial. Then the bundle 2-gerbe
    product gives rise to the following stable isomorphism of bundle gerbes with connection and
    curving:
    \[
    \cQ_{ij} \otimes  \cQ_{jk} \to \cQ_{ik}\otimes \delta(\cG_{ijk})
    \]
    for a bundle gerbe $\cG_{ijk}$ with connection and curving over
    $U_{ijk}=U_i\cap U_j\cap U_k$, hence the curving $B_{ij} = (s_i, s_j)^*B$
    satisfies
    \[
    dB_{ij} = C_i-C_j, \qquad B_{ij} +B_{jk} + B_{ki} = d A_{ijk}
    \]
    for a connection 1-form $A_{ijk}$ on $\cG_{ijk}$. Moreover, the associator $\phi$
    defines a $U(1)$-valued function $g_{ijkl}$ over $U_{ijkl}$ such that
    $g^{\ }_{ijkl}$ satisfies the \v{C}ech 3-cocycle condition
    \[
    g^{\ }_{ijkl}g_{ijkm}^{-1}g^{\ }_{ijlm} g_{iklm}^{-1}g^{\ }_{jklm} =1
    \]
    and
    \[
    A_{ijk}-A_{ijl}+A_{ikl}-A_{jkl}= g_{ijkl}^{-1}d g_{ijkl}.
    \]

\begin{definition}
    A bundle 2-gerbe $(\cQ, Y; X, M)$  is called trivial if $(\cQ, Y; X^{[2]})$
    is isomorphic in $\mathbf{BGrb}_{X^{[2]}}$ to 
    \[
   \delta (\cG) =  \pi_2^*(\cG^*)\otimes \pi^*_1(\cG)
    \]
 for a
 bundle gerbe $\cG$ over $X$ together with compatible conditions 
  on  bundle 2-gerbe
   products  in   $\mathbf{BGrb}_{X^{[3]}}$ and 
 associator natural transformations in $\mathbf{BGrb}_{X^{[4]}}$
  (see \cite{Ste} for more details).  
 A bundle 2-gerbe $(\cQ_1, Y_1; X, M)$ is called stably isomorphic
 to a bundle 2-gerbe $(\cQ_2, Y_2; X, M)$ if and only if
 $\cQ_1$ is isomorphic to $\cQ_2\otimes\delta (\cG)$
    for a bundle gerbe $\cG$ over $X$ together with extra conditions 
    involving the associator natural isomorphisms for $\cQ_1$ and 
    $\cQ_2$.  
    \end{definition}

    \begin{lemma}\label{lemma4}
  Let $(\cP,X; Y,M)$ be a bundle 2-gerbe with connection and curving. Suppose there exists a stable
  isomorphism of bundle gerbes $(\cP,X) \cong (\cQ,Z)$ over $Y^{[2]}$. Then
  there exists a bundle 2-gerbe structure $(\cQ,Z,Y,M)$ with induced connection and curving
  which  has the same Deligne class in $H^3(M,\cD^3)$ as the original bundle
  2-gerbe $(\cP,X; Y,M)$.
  \end{lemma}
  \begin{proof}
  First we must show that $(\cQ,Z; Y,M)$ admits a bundle 2-gerbe product. We use
  the bundle 2-gerbe product on $(\cP,X; Y,M)$ to define it. Recall that this
  product is a stable isomorphism of bundle gerbes over $Y^{[3]}$,
  \[
  \pi_1^*\cP \otimes \pi_3^* \cP \cong \pi_2^* \cP
  \]
  It is convenient here to realise the stable isomorphism as a trivial
  bundle gerbe by expressing the product as a bundle gerbe isomorphism
  \[
  \pi_1^*\cP \otimes \pi_2^* \cP^* \otimes \pi_3^* \cP
  = \delta (J)
  \]
  Similarly we have an isomorphism
  \[
  \cP = \cQ \otimes \delta(L)
  \]
  representing the stable isomorphism of bundle gerbes over $Y^{[2]}$. Thus we
  have
  \[
  \pi_1^*(\cQ \otimes \delta (L)) \otimes \pi_2^* (\cQ \otimes \delta (L))^*
  \otimes \pi_3^* (\cQ \otimes \delta(L))
  = \delta (J)
  \]
  and so
  \[
  \pi_1^*\cQ \otimes \pi_2^* \cQ^* \otimes \pi_3^* \cQ
  = \delta (J) \otimes \pi_1^* \delta(L)^* \otimes \pi_2^*
  \delta(L) \otimes \pi_3^* \delta(L)^*
  \]
  where we use the fact that $\pi_i^*\delta(L) \otimes \pi_i^*\delta(L)^*$
  is canonically trivial. Since the pullback of a trivial bundle gerbe must
  itself be trivial and a tensor product of trivial bundle gerbes is trivial
  then the right hand side is trivial and thus we potentially have a
  bundle 2-gerbe product for $\cQ$.   
  To confirm that it does define a  bundle 2-gerbe product
   we must check the associativity conditions.
   
   Note that it is helpful now to look at diagram (\ref{2-arrow}) and
   diagram (\ref{2-arrow:commute}) to understand the following arguments.
  Recall that there is a bundle
  called the associator bundle on $X^{[4]}$ which is the obstruction to the
  bundle 2-gerbe product being associative. It can be defined by considering
 the product
  \[
  \pi_1^{-1} \delta (J) \otimes \pi_2^{-1} \delta (J)^*
  \otimes \pi_3^{-1} \delta (J) \otimes \pi_4^{-1} \delta (J)^*
  \]
  where $\pi_i : Y^{[4]} \rightarrow Y^{[3]}$ are the face maps in the simplicial
  complex. This product defines the associator  line bundle on
  $Y^{[4]}$ (Cf. Remark \ref{some:mercy} and Diagram (\ref{2-arrow})). 
  Changing to the trivialisation
  representing the bundle gerbe product for $\cQ$,
   we find that the extra terms  involving $\delta(J)$ all cancel
   (in the sense of having   canonical trivialisations), hence the associator
   line bundles for $\cP$ and $\cQ$ are the same, so
  the bundle 2-gerbe product for $\cQ$ is well defined.
  With the induced connection and curving, it is straight forward to show that the
  Deligne class is cohomologous to the the
  Deligne class for $(\cP,X; Y,M)$.
  \end{proof}

 From Lemma \ref{lemma4}, we know that two stably isomorphic bundle 2-gerbes
 with connection and curving have the same Deligne class. Given a representative of
 a Deligne class as in (\ref{local:b2g:d-class}), we can construct a local bundle 2-gerbe
 with connection and curving over $M$ as in \cite{Ste} and \cite{Joh}.
 Analogous to the fact that $H^2(M, \cD^2)$
classifies stable equivalence classes
   of bundle gerbes with connection and curving, we have the following
   proposition, whose complete proof can be found in \cite{Joh}.

   \begin{proposition}(Cf. \cite{Joh})
    \label{b2g:class}
     The group of stable isomorphism classes of bundle 2-gerbes with connection
    and curving over $M$ is isomorphic to $H^3(M, \cD^3)$.
    \end{proposition}

  \section{Multiplicative bundle gerbes}\label{section:multiplicative-gerbe}

   The simplicial manifold $BG_{\bullet}$ associated to the classifying space of $G$
  is constructed in \cite{Dup}, where the total space of the universal $G$-bundle
  $EG$   also has   a simplicial manifold structure.
  The simplicial manifold
  \[
  BG_{\bullet} =\{BG_n =  G \times \cdots \times G \ \text{(n copies)} \}
  \]
   (where $n=0, 1, 2,
  \cdots$), is endowed with  face operators $d_i: G^{n+1} \to G^n$,
  ($i= 0, 1, \cdots, n+1$)
  \[
  d_{i}(g_{0},\ldots,g_{n}) = \begin{cases}
                                  (g_{1},\ldots,g_{n}), & i = 0, \\
                                  (g_{1},\ldots,g_{i-1}g_{i},g_{i+1},
                                    \ldots,g_{n}), & 1\leq i\leq n,
  \\
                                   (g_{0},\ldots,g_{n-1}), & i = n+1.
                                  \end{cases}
  \]

    \begin{definition}
    A {\sl multiplicative bundle gerbe over a compact Lie group $G$} is defined to be a
    simplicial bundle gerbe on the simplicial manifold $BG_{\bullet}$
    associated to the classifying space of $G$.
    \end{definition}

    For a compact, simply connected, simple Lie group $G$,  the tautological
    bundle over $G$ associated to any class in $H^3(G, \ZZ)$ is a simplicial bundle
    gerbe as shown in \cite{Ste}, hence, a multiplicative bundle gerbe.

    \begin{proposition}\label{prop: H^4(BG) and simp bgs on BG}
Let $G$ be a compact, connected Lie group.  Then there is
an isomorphism between $H^4(BG;\ZZ)$
and  the space of isomorphism classes of multiplicative bundle gerbes
on  $G$.
\end{proposition}

First of all, it is not very hard to see that $H^4(BG;\ZZ)$
corresponds to isomorphism classes of simplicial bundle
gerbes on the simplicial manifold $EG^{[\bullet]}$.  This is
because a simplicial bundle gerbe on $EG^{[\bullet]}$ is
the same thing as a bundle $2$-gerbe on $EG\to BG$.   Here
we say that two simplicial bundle gerbes $\cG$ and $\cQ$ on
$EG^{[\bullet]}$ are isomorphic if there is a stable isomorphism
$\cG\cong \cQ$ which is compatible with all the multiplicative structures on
$\cG$ and $\cQ$.  On a general simplicial manifold $X_\bullet$
the notion of isomorphism of simplicial bundle gerbes is
more complicated, involving bundle gerbes on $X_0$, however
because $EG$ is contractible we may use this simpler
notion of isomorphism without any loss of generality.

    Recall from \cite{Bry} and  \cite{BryMcL1} the definition of the simplicial
\v{C}ech cohomology groups $H^*(X_\bullet;\underline{A})$
for a simplicial manifold $X_\bullet$ and some topological abelian group $A$.
 To define these one first  needs the notion
of a covering of the simplicial manifold $X_\bullet$.  By definition
this is a family of covers $\cU^n = \{U^n_\a\}$ of
the
manifolds $X_n$ which are compatible with the face and degeneracy
operators for $X_\bullet$.  Brylinski and McLaughlin in \cite{BryMcL1}
explain how one may inductively construct such a family of
coverings by first starting with an arbitrary cover $\cU^0$ of
$X_0$ and then choosing a common refinement $\cU^1$ of the induced
covers $d_0^{-1}(\cU^0)$ and $d_1^{-1}(\cU^0)$ of $X_1$.  $\cU^1$
then induces three covers $d_0^{-1}(\cU^1)$, $d_1^{-1}(\cU^1)$
and $d_2^{-1}(\cU^1)$ of $X_2$.  One then chooses a common
refinement $\cU^2$ and repeats this process.  In particular, the
covering $\cU^\bullet$ may be chosen so that each $\cU^n$ is a good
cover of $X_n$.

The simplicial \v{C}ech cohomology
$H^*(\cU^{\bullet},\underline{A})$ of $X_\bullet$ for the
covering
$\cU^\bullet$ is by definition the cohomology of the double complex
$C^p(\cU^q,\underline{A})$ where $C^*(\cU^p,
\underline{A})$ is the \v{C}ech complex for the covering $\cU^p$
of the manifold $X_p$.  The differential for the complex
$C^p(\cU^*,\underline{A})$ is induced in the usual way
from the face operators $d_i$ on $X_\bullet$.  The groups
$H^*(X_\bullet,\underline{A})$ are then defined by taking a direct limit of the
coverings $\cU^\bullet$.  If the covering $\cU^\bullet$ is good
in the sense that each $\cU^n$ is a good cover then
$$
H^*(X_\bullet;\underline{A}) \cong
H^*(\cU^\bullet;\underline{A}).$$
Note that there is a spectral sequence converging to a graded
quotient of $H^*(X_\bullet;\underline{A})$ with
$$
E_1^{p,q} = H^p(X_q;\underline{A}).
$$
The following proposition is a straightforward extension
of Theorem 5.7 in part I of \cite{BryMcL1}  to the
language of bundle gerbes.

\begin{proposition}\label{BM}
Let $X_\bullet$ be a simplicial manifold.  Then we have that
isomorphism classes of simplicial bundle gerbes on $X_\bullet$
are classified by the simplicial \v{C}ech cohomology group
$H^3(X_{\bullet \geq 1};\underline{U(1)})$.
Here $X_{\bullet \geq 1}$ denotes the truncation of $X_\bullet$
through degrees $\geq 1$.
\end{proposition}

We sketch a proof of this Proposition below.  Let us first be clear
about what we mean by the group $H^3(X_{\bullet \geq 1};
\underline{\CC}^*)$.  By this we mean that if $\cU^\bullet$ is a
good covering of $X_\bullet$, so that $H^*(X_\bullet,
\underline{U(1)}) = H^*(\cU^\bullet,
\underline{U(1)})$, then $H^*(X_{\bullet \geq 1},
\underline{U(1)})$ is the cohomology of the double complex
$C^p(\cU^q;\underline{U(1)})$ with $q$ in degrees $\geq 1$.
Also, by an isomorphism of simplicial bundle gerbes on $X_\bullet$ we
mean a stable isomorphism $\cG\cong \cQ$ compatible with all the product
structures
on $\cG$ and $\cQ$.  We do not require that $\cG\cong \cQ\otimes \d(T)$
for some bundle gerbe $T$ on $X_0$.  As noted above this is a restrictive
definition but is sufficient for the cases we are interested in, such as the
contractible space $EG$.

Given a simplicial bundle gerbe $\cG$ on $X_\bullet$, we associate to
$\cG$ a simplicial \v{C}ech cohomology class in
$H^3(X_{\bullet \geq 1},\underline{U(1)})$ as follows.
For the good covering $\cU^\bullet$ on $X_\bullet$, let $\underline{g} =
(g_{\a\b\c})$ be a \v{C}ech
cocycle representative for the Dixmier-Douady class of $\cG$.  Then it
is easy to see that the $2$-cocycle
\[
d_0^*(\underline{g})\,
d_1^*(\underline{g}^{-1})
d_2^*(\underline{g}) = \d(\underline{h})
\]
for some $1$-cochain
$\underline{h} = (h_{\a\b})$ on the covering $\cU^2$.
Then the simplicial \v{C}ech  $1$-cochain
\[
d_0^*(\underline{h})d_1^*(\underline{h}^{-1})d_2^*(\underline{h})
d_3^*(\underline{h}^{-1})
\]
 on the cover $\cU^3$ is a $1$-cocycle: it
is a \v{C}ech representative for the first Chern class of the associator line
bundle on $X_3$ (the line bundle induced by the associator).
Consequently, we must have
\[
d_0^*(\underline{h})\,d_1^*(\underline{h}^{-1})
d_2^*(\underline{h})\,d_3^*(\underline{h}^{-1}) = \d(\underline{k})
\]
 for some simplicial \v{C}ech  $0$-cochain $\underline{k} = (k_\a)$ on $\cU^3$.  The
cocycle condition for the associator section shows that we must have
\[
d_0^*(\underline{k})\, d_1^*(\underline{k}^{-1}) d_2^*(\underline{k})\,
d_3^*(\underline{k}^{-1})d_4^*(\underline{k}) = 1
\]
 on $\cU^4$.
The triple $(\underline{g},\underline{h},\underline{k})$ is a simplicial \v{C}ech  cocycle
in the truncated double complex $C^p(\cU^{q\geq 1},\underline{U(1)})$
representing a class in $H^3(X_{\bullet \geq
1},\underline{U(1)})$.
Conversely, such a simplicial \v{C}ech  cocycle $(\underline{d},\underline{h},\underline{k})$
in the truncated double complex $C^p(\cU^{q\geq 1},\underline{U(1)})$
determines a unique isomorphism class of simplicial bundle gerbes on
$X_\bullet$.  If the cocycle
$(\underline{d},\underline{h},\underline{k})$
associated to a simplicial bundle gerbe $\cG$ is trivial in $H^3(X_{\bullet \geq
1},\underline{U(1)})$, then there is a
stable isomorphism between $\cG$ and the trivial simplicial bundle gerbe,
consisting of a trivial bundle gerbe on $X_1$ equipped with the trivial
product structures.

As an immediate consequence of  Proposition \ref{BM}, we see that
isomorphism classes of simplicial bundle gerbes on $BG_\bullet$
are classified by the simplicial \v{C}ech cohomology group
$H^3(BG_\bullet,\underline{U(1)})$.  Since the simplicial map
$EG^{[\bullet]} \to BG_\bullet$ is a homotopy equivalence in each
degree, we see from the spectral sequence above that it induces
an isomorphism $$H^3(BG_\bullet,\underline{U(1)}) \cong
H^3(EG^{[\bullet]};\underline{U(1)}).$$
As we have already  noted, isomorphism
classes of simplicial bundle gerbes on $EG^{[\bullet]}$ correspond
exactly to
isomorphism classes of bundle $2$-gerbes on $EG\to BG$
and hence to $H^4(BG;\ZZ)$.

  Given a principal $G$-bundle $\pi: P \to M$, there exists a natural map
  \[
  \hat{g}:  P^{[2]} \longrightarrow G
  \]
  given by $p_1\cdot \hat{g} (p_1, p_2) = p_2$ such that $P^{[2]} \to P \times G$
 sending $(p_1, p_2)$ to $ (p_1, \hat{g} (p_1, p_2)$ is a diffeomorphism. 
  
  \begin{lemma}\label{b2g:M}
  Given a principal $G$-bundle $P\rightarrow M$ and a multiplicative
  bundle gerbe $\cG$ over $G$ then there exists a bundle
  2-gerbe over $M$ of the form $(\cQ, X; P, M)$ such that $(\cQ, X, P^{[2]})$
  is the pull-back bundle gerbe $\hat{g}^*\cG$.
  \end{lemma}
  \begin{proof} Note that there exists a diffeomorphism 
  \[
  P \times G^n  \longrightarrow P^{[n+1]}
  \]
  given by 
  \[
  (p, g_1, g_2, \cdots, g_n) \mapsto \bigl(p, p\cdot g_1,  p\cdot g_1g_2,
  \cdots, p\cdot (g_1g_2\cdots g_n)\bigr).
  \]
  The inverse of this map, composing with the projection to $G^n$,
   defines to a simplicial map
  \[\hat{g}_\bullet: 
  P_\bullet =\{P_n = P^{[n+1]}\} \longrightarrow BG_{\bullet}=\{ BG_n = G^n\},
  \]
  which can be used to pull back the simplicial bundle gerbe
  over $BG_{\bullet}$ corresponding to $\cG$ to a simplicial bundle gerbe 
  $\cQ$ over $P_\bullet$. 
  This defines a bundle 2-gerbe $(\cQ, Y; P, M)$ over $M$ of the required form. 
  \end{proof}
  
Given a multiplicative  bundle gerbe $\cG$ over $G$, applying the above Lemma \ref{b2g:M}
 to the universal bundle $\pi: EG\to BG$, we obtain
 a bundle 2-gerbe over $BG$ of the form $(\hat{g}^*\cG; EG,BG)$
 with  $\hat{g}^*\cG$ is a bundle gerbe over $EG^{[2]}$ obtained by
 the pull-back of $\cG$ via $\hat{g}: EG^{[2]} \to G$. Conversely,
 we will show that every bundle 2-gerbe over $BG$ is stably isomorphic a
 bundle 2-gerbe of this form.
 
\begin{lemma}\label{lemma2}
  Every bundle 2-gerbe on $BG$ is stably isomorphic to a bundle 2-gerbe
  of the form $(\cQ,X; EG,BG)$, where $(\cQ,X)$ is a bundle
  gerbe over $EG^{[2]}$.
  \end{lemma}
  \begin{proof}
  We use the classifying theory of bundle 2-gerbes. It is well known that
  $H^4(M,\ZZ) \cong [M;K(\ZZ, 4)]$. We use the iterated classifying space
  $BBBU(1)$ (or $B^3U(1)$) as a model for $K(\ZZ,4)$ with a differential space
  structure constructed in \cite{Gaj}. In  \cite{Gaj} Theorem H , it is shown 
  that for a smooth manifold $M$, the group $H^4(M, \ZZ)$ is isomorphic to the group
  of isomorphism classes of smooth principal $B^2U(1)$-bundles over $M$.
  
  We can    transgress the degree $4$ class in $H^4(B^3 U(1), \ZZ)$
  to get a degree 3 class in $H^3(B^2U(1), \ZZ)$ which determines a 
  multiplicative bundle gerbe over $B^2U(1)$. Then we  apply the
  canonical map $\hat{g}: {EB^2 U(1)}^{[2]}\to B^2 U(1)$ to pull-back the
  corresponding multiplicative bundle gerbe over $B^2U(1)$ 
  to ${EB^2 U(1)}^{[2]}$. This gives rise to  the universal
  bundle 2-gerbe $\tilde\cQ$ over $B^3 U(1)$.
     The classifying bundle   2-gerbe then has the form
  \[
  \begin{array}{ccc}
  \tilde{\cQ }& & \\
  \Downarrow & & \\
  {EB^2 U(1)}^{[2]} & \rightrightarrows & EB^2U(1) \\
  & & \downarrow \\
  & & B^3 U(1)
  \end{array}
  \]
  As $B^3 U(1)$ is 3-connected, and $H^4(B^3 U(1), \ZZ) \cong \ZZ$,
   the tautological bundle 2-gerbe developed in \cite{CMW} can be adapted 
   to give another construction of
     such a classifying bundle  2-gerbe over $B^3 U(1)$ associated to
  any integral class in $H^4(B^3 U(1), \ZZ)$.

  Any bundle 2-gerbe on $BG$ is defined by pulling back the universal
  bundle 2-gerbe by a classifying map $\psi: BG \rightarrow B^3U(1)$. Consider
  the map $\pi \circ \psi : EG \rightarrow B^3U(1)$ where $\pi$ is the projection
  in the universal $G$-bundle. By using the homotopy lifting property and the
  contractibility of $EG$ we can  always find a lift $\hat{\psi}$,
  \[
  \begin{array}{ccc}
  EG & \stackrel{\hat{\psi}}{\rightarrow} & EB^2U(1) \\
  \downarrow & & \downarrow \\
  BG & \stackrel{\psi}{\rightarrow} & B^3U(1)
  \end{array}
  \]
  Thus we can pull back the universal bundle 2-gerbe to get a bundle
  2-gerbe of the form $(\cQ,X; EG,BG)$,
  with $(\cQ,X)$  a bundle   gerbe over $EG^{[2]}$.
  \end{proof}

  \begin{lemma}\label{lemma3}
  Any bundle gerbe over $EG^{[2]}$ is stably isomorphic to $\hat{g}^* \cG$ where
  $\hat{g}: EG^{[2]} \rightarrow G$ is the map satisfying $e_2 = e_1
  \hat{g}(e_1,e_2)$ for $(e_1,e_2) \in EG^{[2]}$ and $\cG$ is some bundle gerbe
  over $G$.
  \end{lemma}
  \begin{proof}
  We may identify $EG^{[2]}$ with $EG \times G$ via the map
  $(e_1,e_2) \mapsto (e_1,\hat{g}(e_1,e_2))$. Thus stable isomorphism
  classes of bundle gerbes over $EG^{[2]}$ are classified by
  $H^3(EG \times G,\ZZ)$ which, since $EG$ is contractible, 
  is equal to $H^3(G;\ZZ)$. Under the identification we see that the
  Dixmier-Douady class of the bundle gerbe on $EG^{[2]}$ must be obtained
  from a class in $H^3(G,\ZZ)$ by the map $\hat{g}$.
  \end{proof}

  \begin{proposition}\label{bgprop}
   Every bundle 2-gerbe on $BG$ is stably isomorphic to a bundle 2-gerbe
   of the form $(\hat{g}^*\cG, X; EG, BG)$ for a multiplicative
   bundle gerbe $\cG$ over $G$.
  \end{proposition}
  \begin{proof}
  We start with any bundle 2-gerbe on $BG$. By Lemma \ref{lemma2} we may
  assume without loss of generality that it is of the form $(\cQ,X; EG,BG)$.
  Next we use lemma \ref{lemma3} to replace the bundle gerbe $(\cQ,X,EG^{[2]})$
  with $(\hat{g}^*\cG,\hat{g}^*Y,EG^{[2]})$ where
  $(\cG, Y, G) $ is now a bundle gerbe over
  $G$. Note that $(\cQ,X; EG^{[2]})$
  and $(\hat{g}^*\cG,\hat{g}^*Y,EG^{[2]})$ is stably isomorphic.
  Using Lemma \ref{lemma4}, we know that there exists a bundle 2-gerbe structure on
  $(\hat{g}^*\cG,\hat{g}^*Y;,EG, BG)$ which
  does not change  the stable isomorphism class of the
  original bundle 2-gerbe on $BG$.
  For the bundle 2-gerbe $(\hat{g}^*\cG,\hat{g}^*Y; EG, BG)$ to be defined,
  $\hat{g}^*\cG$  must be a simplicial bundle gerbe over
   $EG_\bullet = \{ EG_n = EG^{[n+1]}\}$.
  Fix a point $p \in EG$, then $\hat{p}: G^n \to EG^{[n+1]}$ with
   \[
  \hat{p} ( g_1, g_2, \cdots, g_n)= \bigl(p, p\cdot g_1,  p\cdot g_1g_2,
  \cdots, p\cdot (g_1g_2\cdots g_n)\bigr),
  \]
  is a simplicial map between $BG_\bullet= \{BG_n = G^n\}$ and 
  $EG_\bullet$. Then $\cG \cong \hat{p}^* \circ \hat{g}^* \cG$
  implies that  $\cG$ is a simplicial bundle gerbe over $BG_\bullet$. Hence,
 $\cG$ is a multiplicative bundle gerbe over $G$.
  \end{proof}

  \begin{theorem}\label{trans=multi:chern}
 Given a bundle gerbe $\cG$ over $G$, $\cG$ is multiplicative if and only if 
 its Dixmier-Douady class is transgressive, i.e., in the image of the transgression
 map   $\tau: H^4(BG,\ZZ) \rightarrow H^3(G,\ZZ)$.
  \end{theorem}
  \begin{proof} Given a multiplicative  bundle gerbe $\cG$
  over $G$, applying  Lemma \ref{b2g:M}
 to the universal bundle $\pi: EG\to BG$, we obtain
 a bundle 2-gerbe $\cQ$  over $BG$ of the form $(\hat{g}^*\cG; EG,BG)$
 with  $\hat{g}^*\cG$ is a bundle gerbe over $EG^{[2]}$ obtained by
 the pull-back of $\cG$ via $\hat{g}: EG^{[2]} \to G$. The stable
 isomorphism class of $(\hat{g}^*\cG; EG,BG)$ defines a class $\phi\in H^4(BG, \ZZ)$.
 We now show that the transgression of $\phi$ under the map 
 $\tau: H^4(BG,\ZZ) \rightarrow H^3(G,\ZZ)$ is the Dixmier-Douady class
 of $\cG$, hence transgressive.
 
 The class $\phi$ defines a homotopy class in $[BG, K(\ZZ, 4)]$ such that 
 the bundle 2-gerbe $\cQ$ is stably isomorphic to a bundle 2-gerbe obtained
  by pulling back the universal
  bundle 2-gerbe $\tilde \cQ$ over $B^3U(1)=K(\ZZ, 4)$ via a classifying map 
  $\psi: BG \rightarrow K(\ZZ, 4)$ of $\phi$
  and the commutative diagram
  \[
  \begin{array}{ccc}
  EG & \stackrel{\hat{\psi}}{\rightarrow} & EK(\ZZ, 3) \\
  \downarrow & & \downarrow \\
  BG & \stackrel{\psi}{\rightarrow} & K(\ZZ, 4).
  \end{array}
  \]
  This commutative diagram gives rise to
  a homotopy class of maps $\hat{\psi}: G \to  K(\ZZ, 3)$
  which determines a class 
  \[
  c_\psi \in H^3(G, \ZZ) \cong [G, K(\ZZ, 3)].
  \]

  Using the long exact sequence for homotopy groups
  for the $K(\ZZ, 3)$-fibration  $EK(\ZZ, 3)\to BK(\ZZ, 3) = K(\ZZ, 4)$ and the
  Hurewicz theorem for 3-connected spaces, we know that
  the transgression map
  $
  H^4(BK(\ZZ, 3), \ZZ) \to H^3(K(\ZZ, 3), \ZZ) $
  sends the generator of $ H^4(BK(\ZZ, 3), \ZZ) $
  to the generator of $H^3(K(\ZZ, 3), \ZZ)$. Therefore, we 
  obtain that
  \[
  c_\psi = \tau (\phi),
 \]
 as $\phi \in H^4(BG, \ZZ) \cong [BG, K(\ZZ, 4)]$ is defined by 
 pulling back the generator of  $ H^4(BK(\ZZ, 3), \ZZ) $
 via the  classifying map $\psi$. 
 
 As  bundle gerbes over $EG^{[2]}$, $\cQ=\hat{g}^* \cG $ is stably
 isomorphic the pull-back bundle gerbe $(\hat{\psi}^{[2]})^*\tilde{\cQ}$
 via the map $\hat{\psi}^{[2]}:
 EG^{[2]}\to EK(\ZZ, 3)^{[2]}$. This implies that
 the Dixmier-Douady class of $\hat{g}^*\cG$ is given by the homotopy
 class of the map
 \[
 \xymatrix{
 EG^{[2]}\ar[r]^{\hat{\psi}^{[2]}}& EK(\ZZ, 3)^{[2]}\ar[r] &K(\ZZ, 3).}
 \]
 As $\hat{g}^*$ induces an isomorphism 
 $H^3(G, \ZZ) \to H^3(EG^{[2]}, \ZZ)$, putting all these together, we know that
 the Dixmier-Douady class of $\cG$ is given by  $c_\psi=\tau(\phi)$, hence, 
 transgressive.
 
  Conversely, suppose $\cG$ is a bundle gerbe  over $G$ whose
  Dixmier-Douady class is a  transgressive  class $\tau (c) \in H^3(G, \ZZ)$ for
  a class   $c \in H^4(BG,\ZZ)$.
  By Proposition \ref{bgprop} we may, without loss of generality, realise
   $c$  by a bundle 2-gerbe $(\cQ_c, X; EG, BG)$
   over $BG$ from  a multiplicative bundle 
  gerbe $\cG_c$ over $G$ such that $\cQ_c = \hat{g}^* \cG_c$. 
  The above argument shows that the Dixmier-Douady class
   of $\cG_c$ is also given by $\tau (c)$. Therefore,
   the bundle gerbe $\cG$ is stably isomorphic to the multiplicative
    bundle  gerbe $\cG_c$. Then  the multiplicative
    structure on $\cG_c$
    \[
    m_{_c}:  d_0^* \cG_c \otimes d_2^* \cG_c \to d_1^* \cG_c
    \]
         induces a multiplicative structure on $\cG$
	 \[
	m:   d_0^* \cG  \otimes d_2^* \cG  \to d_1^* \cG
	  \]
     as a stable isomorphism in the bi-category $\mathbf{BGrb}_{G^2}$.
      The  associator for $(\cG, m_\cG)$ (see (\ref{2-arrow}))
      \[
      \phi:   d_2^*m\circ (d^*_0m\otimes Id) 
      \to d_1^*m \circ (Id \otimes d^*_3 m)
      \]
      is also induced by the corresponding associator $\phi_c$
       for $(\cG_c, m_{c})$ in the bi-category $\mathbf{BGrb}_{G^3}$.
      The coherent condition for $\phi$ (see (\ref{2-arrow:commute}))
      follows from  the corresponding 
   coherent condition for $\phi_c$ in the bi-category $\mathbf{BGrb}_{G^4}$.
   Hence, the bundle gerbe $\cG$ over $G$,  whose
  Dixmier-Douady class is  transgressive, is multiplicative.   
  \end{proof}

   Note that there exists a simpler proof of Theorem \ref{trans=multi:chern}
   by applying Proposition    \ref{prop: H^4(BG) and simp bgs on BG} and 
   an  observation that
   $H^4(BG, \ZZ)$ classifies the isomorphism classes of simplicial bundle
   gerbe over $BG_\bullet$. The constructive proof for Theorem \ref{trans=multi:chern}
   is given so that the proof can be adopted to establish the following
    theorem, which is the refinement of
    Theorem \ref{trans=multi:chern} for the correspondence
    from three dimensional Chern-Simons gauge theories with gauge
    group $G$ to Wess-Zumino-Witten models on the group manifold $G$.
    The proof will be postponed until after
    we discuss  Chern-Simons bundle 2-gerbes.

   \begin{theorem}\label{CS:image}
    Denote by
    $
    \Psi: CS(G) \longrightarrow WZW (G) \cong H^2(G,\cD^2)
    $
    the correspondence map  and let $\cG$ be a bundle
    gerbe over $G$ with connection and curving, whose Deligne class $d(\cG)$
    is in $H^2(G, \cD^2)$. Then
    $d(\cG) \in Im (\Psi)$, if and only if $\cG$ is multiplicative.
    \end{theorem}

\begin{remark} Note that  the transgression map $H^4(BSO(3), \ZZ)
\to H^3(SO(3), \ZZ)$ sends the generator of $H^4(BSO(3), \ZZ)$ 
to twice of the generator of  $H^3(SO(3), \ZZ)$. 
Hence, only those bundle gerbes $\cG_{2k}$  over $SO(3)$, with even
Dixmier-Douady classes in $ H^3(SO(3), \ZZ)\cong \ZZ$, are multiplicative. 
\end{remark}

  \section{The Chern-Simons  bundle 2-gerbe}
  \label{CS2gerbe}

  In this Section, we will construct a Chern-Simons bundle 2-gerbe
  over $BG$ such that the proof of the main theorem (Theorem \ref{CS:image}) follows.

  Given a principal $G$-bundle $P$ with a connection $A$ over $M$,
a Chern-Simons gauge
  theory $c\in CS(G)$ defines a degree 3 Deligne  characteristic class
  \[
  c(P, A) \in H^3(M, \cD^3).
  \]
  We will construct a bundle 2-gerbe with connection and curving
  over $M$ corresponding to the Deligne class $c(P, A)$.
We shall call this a {\it Chern-Simons bundle
  2-gerbe}. We define it in terms of a universal Deligne characteristic
class represented
  by the {\it universal Chern-Simons bundle 2-gerbe}.

 In order to use differential forms on $BG$ we need
 to fix a smooth infinite dimensional model of
    $BG$ by embedding $G$ into $U(N)$ and letting $EG$ be the Stiefel manifold of
    $N$ orthonormal vectors in a separable complex Hilbert space. Alternatively,
     we could choose
 a smooth finite dimensional $n$-connected 
 approximation $B_n$  of the classifying space
 $BG$ as in \cite{NarRam}. Given
  a principal $G$-bundle $P$ with a connection $A$ over $M$ and an integer
  $n>max\{5, dim M\}$,  there is a choice of $n$-connected finite dimensional
   principal  $G$-bundle $E_n$ with a connection $\bar A$ such that
  $(E_n, \bar A)$ is a classifying space of $(P, A)$.
  In particular $H^k(B, \ZZ) \cong H^k(BG, \ZZ)$ for $k\leq n$.
  For convenience, we suppress this latter detail and
  work directly on the  infinite dimensional smooth model of $EG\to BG$.

   Note that there exists a connection
  $\AA$ on the universal bundle $EG$ over $BG$ and a classifying map
  $f$ for $(P, A)$ such that
  \[
  (P, A) = f^* (EG, \AA),
  \]
  which implies $c(P, A)= f^* c(EG, \AA)$ with $c(EG, \AA)\in H^3(BG, \cD^3)$
  is the Deligne characteristic class for $(EG, \AA)$. From the commutative diagram:
  \[
  \begin{array}{ccc}
  H^3(BG, \cD^3)&\longrightarrow& H^2(G, \cD^2)\\[2mm]
\downarrow^c && \downarrow^c \\[2mm]
H^4(BG, \ZZ)& \to & H^3(G, \ZZ),
\end{array}
\]
where the vertical maps are the characteristic class maps on Deligne cohomology groups,
we see that Theorem \ref{CS:image} refines the result in Theorem \ref{trans=multi:chern}.

Given a class $\phi \in H^4(BG, \ZZ)$, denote by $\Phi$ the corresponding
$G$-invariant polynomial on the Lie algebra of $G$.  Then  associated to
a connection $A$ on a principal $G$-bundle (with curvature $F_A$)
there is a closed 4-form
\[
\Phi (\disp{\frac{i}{2\pi}} F_A) \in \Omega^4_{cl, 0}(M)
\]
with integer periods. In fact, $\Phi (\disp{\frac{i}{2\pi}} F_A)
= f^*\Phi (\disp{\frac{i}{2\pi}} F_\AA)$, and
\[
\bigl(\phi, \Phi (\disp{\frac{i}{2\pi}} F_\AA) \bigr)\in H^4(BG,\ZZ)
\times \Omega^4_{cl, 0}(BG)
\]
determines a unique Deligne class
\[
c_\phi(EG, \AA) \in H^3(BG, \cD^3).
\]
Hence $c_\phi (P, A) =f^*(c_\phi(EG, \AA)) \in H^3(M, \cD^3)$. So we can
say that a class $\phi \in H^4(BG, \ZZ)$ defines a {\it
canonical Chern-Simons
gauge theory} $c_\phi$ with gauge group $G$.

   From the exact sequence for $H^3(BG, \cD^3)$,
\[
0 \rightarrow
\Omega^3(BG)/ \Omega^3_{cl, 0}(BG) \rightarrow H^3(BG, \cD^3) \rightarrow H^{4}(BG, \ZZ)
  \rightarrow 0,
  \]
 any Chern-Simons gauge theory with the same characteristic  class $\phi$ in
 $H^{4}(BG, \ZZ)$ differs from $c_\phi$ by a Deligne class
 $[1, 0, 0, C]$ for a 3-form $C$ on $BG$. Note that $[1, 0, 0, C]$
 defines a trivial bundle 2-gerbe over $BG$, so in this section, we only construct
 the universal Chern-Simons bundle 2-gerbe over $BG$ corresponding to the
 canonical Chern-Simons gauge theory $c_\phi$ with gauge group $G$.
 Notice that for two different connections
 $\AA_1$ and $\AA_2$ on the universal
 bundle $EG$,
 \[
 c_\phi(EG, \AA_1)-c_\phi(EG, \AA_2) = [(1, 0, 0, CS_\phi(\AA_1, \AA_2))],
 \]
 where $CS_\phi(\AA_1, \AA_2)$ is the Chern-Simons form on $BG$ associated to a pair of connections
 $\AA_1$ and $\AA_2$ on $EG$ (Cf. \cite{Chern}).

 \begin{remark}\label{S:P-A}
 Recall \cite{CheSim} that associated with each principal
 $G$-bundle $P$ with connection $A$  Cheeger and Simons constructed a
  differential character $$S_{\Phi,\phi} (P, A) \in \check{H}^3(M,U(1)), $$ where
  $\Phi \in I^2(G)$ is a $G$-invariant polynomial on its Lie algebra
   and $\phi \in H^4(BG,\ZZ)$ is  a characteristic class corresponding to $\Phi$  under
  the Chern-Weil homomorphism. This differential
  character is uniquely defined when it satisfies the following:
  \begin{enumerate}
  \item The image of $S_{\Phi,\phi}(P, A)$ under the curvature map
$$\check{H}^3(M,U(1))
  \rightarrow \Omega_{cl,0}^4(M,\RR)$$ is $\Phi(\disp{\frac{i}{2\pi}}F_A)$ where $F_A$ is the curvature
  form of $A$.
  \item The image of $S_{\Phi,\phi}(P, A)$ under the characteristic class map $\check{H}^3(M,U(1))
  \rightarrow H^4(M,\ZZ)$ is $\phi(P)$, the characteristic class of $P$ associated to $\phi$.
  \item The assignment of $S_{\Phi,\phi}(P, A)$ to $(P, A)$ is natural with respect
  to morphisms of principal $G$-bundles with connection.
  \end{enumerate}
  Since differential characters and bundle 2-gerbes with connection and curving
   on $M$ are both classified   by the Deligne cohomology group $H^3(M,\cD^3)$,
    our Chern-Simons bundle 2-gerbe is a bundle gerbe version
  of the Cheeger-Simons invariant described above.
  \end{remark}

  Given a connection $\AA$ on $\pi: EG\to BG$, the canonical Chern-Simons
gauge theory $c_\phi$
   associated to
  a class $\phi \in H^4(BG, \ZZ)$
defines a universal Deligne characteristic class
  $c_\phi(EG, \AA)$  from the pair  $(\phi, \Phi(\disp{\frac{i}{2\pi}}F_\AA))$. 
  Then there is a Chern-Simons form
  $CS_\phi (\AA)$ associated
  to $c_\phi$ and $(EG, \AA)$, satisfying
  \ba\label{CS:form}
  d CS_\phi (\AA) = \pi^* \bigl(\Phi(\disp{\frac{i}{2\pi}}F_\AA)\bigr),
  \na
 and  the restriction of $CS_\phi (\AA)$
  to a fiber of $EG\to BG$ determines a left-invariant closed 3-form $\omega_\phi$ on $G$.

  This universal Chern-Simons form $CS_\phi (\AA)$ can be
  constructed as in \cite{Fre} via the pull-back
  principal $G$-bundle $\pi^* EG \to EG$, which admits a section $e\mapsto (e, e)$.
  This trivialisation defines a trivial connection $\AA_0$ on $\pi^* EG$. Then
  $\pi^*\AA$ and $\AA_0$ defines a path of connections
  \[
  \AA_t = t\pi^*\AA + (1-t) \AA_0
  \]
  for $0\leq t \leq 1$, which can be  thought as a connection on
  $[0, 1]\times  \pi^* EG \to [0, 1]\times  EG$. Define
  \ba\label{CS:form1}
  CS_\phi(\AA) =\disp{\int_{[0, 1]} \Phi(\frac{i}{2\pi}}F_{\AA_t}).
  \na
  Then the relation (\ref{CS:form}) follows from Stokes' theorem for the
  projection $[0, 1]\times  EG \to EG$ and $\Phi(\frac{i}{2\pi}F_{\AA_0}) =0$
  for  $\AA_0$  a trivial connection. 
  
  \begin{remark}The corresponding left-invariant
   closed 3-form $\omega_\phi$ on $G$ is an integer multiple of the standard
   3-form $<\theta, [\theta, \theta]>$ where $\theta$ is the left-invariant
   Maurer-Cartan form on $G$ and $<\ , \ >$ is the symmetric bilinear form on the
   Lie algebra of $G$ defined by $\Phi\in I^2(G)$.
   \end{remark}

    Given the fibration $EG\to BG$, we introduce the natural map
$\hat{g}: EG^{[2]}\to G$
  defined by $e_2 = e_1\cdot \hat{g} (e_1, e_2)$ where $(e_1, e_2)
  \in EG^{[2]}$.

 \begin{definition}\label{universal:cs-b2g}
  The
  {\it universal Chern-Simons bundle 2-gerbe} $\cQ_\phi$
associated to
$\phi \in H^4(BG, \ZZ)$ is a bundle 2-gerbe $(\cQ_\phi, EG^{[2]}; EG, BG)$
 illustrated by the following diagram:
  \[
  \xymatrix{
  & \cG \ar@2{->}[d]\\
 \cQ_\phi \ar@2{->}[d]&(G, \omega_\phi) \\
  EG^{[2]}\ar[ur]^{\hat{g}}\ar@< 2pt>[r]^{\pi_1} \ar@< -2pt>[r]_{\pi_2}
  & (EG,CS_\phi (\AA) ) \ar[d]_{\pi} \\
   &(BG, \Phi(\disp{\frac{i}{2\pi}}F_\AA))
  }
  \]
where  the bundle gerbe $\cQ_\phi$ over $EG^{[2]}$ is obtained
from the pull-back of  a multiplicative bundle
  gerbe $\cG$ over $G$ associated to $\tau (\phi)\in H^3(G, \ZZ)$.
  Given   a connection $\AA$ on $EG$, the  bundle gerbe $\cQ_\phi$ over $EG^{[2]}$ is 
  equipped with a connection whose  bundle gerbe curvature is given by
  $(\pi_2^*-\pi^*_1)CS_\phi (\AA)$, with $CS_\phi (\AA)$
  the Chern-Simons form (\ref{CS:form1}) on $EG$ associated to $\phi$ and $\AA$, and the
  bundle 2-gerbe curvature is given by $\Phi(\disp{\frac{i}{2\pi}}F_\AA)$.
Moreover the bundle gerbe $\cG$ over $G$
 is equipped with a connection and curving with the bundle
  gerbe curvature given by $\omega_\phi$.
  \end{definition}

  \begin{proposition}  \label{cs2bg:exist}
   Given a class $\phi \in H^4(BG, \ZZ)$, there exists a universal
  Chern-Simons bundle 2-gerbe over $BG$ associated to a connection $\AA$ on $EG$.
  \end{proposition}
  \begin{proof}
  From the proof of Lemma \ref{lemma2}, Proposition \ref{bgprop} and Theorem
  \ref{trans=multi:chern},
   we can represent a class $\phi\in H^4(BG, \ZZ)$ by
   a bundle 2-gerbe $(\cQ_\phi, X; EG, BG)$
   over $BG$  which is associated to the universal  $G$-bundle $EG\to BG$ and
to  a multiplicative bundle
  gerbe $\cG$ over $G$, and is such that $\cQ_\phi$
 is stably isomorphic to $\hat{g}^*\cG$ for $\hat{g}: EG^{[2]}\to G$.
  To complete the proof we have to equip  $(\cQ_\phi, X; EG, BG)$ with
a bundle 2-gerbe connection and curving.

  Let $g: P \to G$ be a gauge transformation on $P$.
  Denote by $\AA^{g}$ the connection on $EG$ such that
  \[
  \AA^{g} (e) = \AA(e\cdot g)
  \]
  Then by direct calculation , we know that
  \[
  CS_\phi (\AA^g) = CS_\phi (\AA) + g^*\omega_\phi - d\Phi (\AA, dg\cdot g^{-1}),
  \]
  where $CS_\phi (\AA)$ is the universal Chern-Simons form
  (\ref{CS:form1}) associated to
  $\phi$ and $(EG, \AA)$, and the left-invariant closed 3-form $\omega_\phi$ on $G$ (the
  transgression of $\Phi (\disp{\frac{i}{2\pi}}F_\AA)$). Then under the map
  $\pi_2^*-\pi^*_1$, $CS_\phi (\AA)$ is mapped to a closed 3-form on $EG^{[2]}$ with periods
  in $\ZZ$.
To see this note that $(\pi_2^*-\pi^*_1)CS_\phi (\AA)$ is closed, as
  \[
  d (\pi_2^*-\pi^*_1)CS_\phi (\AA) = (\pi_2^*-\pi^*_1)\circ \pi^* \Phi (\disp{\frac{i}{2\pi}}F_\AA)
  =0,
  \]
  and to confirm that $(\pi_2^*-\pi^*_1)CS_\phi (\AA)$ has its periods in $\ZZ$, we take
  a 3-cycle $\sigma$ in $EG^{[2]}$, and form a 4-cycle in $EG$ given
  by gluing two chains in $EG$ with boundary $\pi_2 (\sigma)$ and $-\pi_1(\sigma)$
  as $EG$ is contractible. Then
  $\Phi (\disp{\frac{i}{2\pi}}F_\AA) \in \Omega^4_{cl, 0}(BG)$ implies that
  $(\pi_2^*-\pi^*_1)CS_\phi (\AA) \in \Omega^3_{cl, 0}(EG^{[2]})$. Actually, we have
  an explicit expression for   $(\pi_2^*-\pi^*_1)CS_\phi (\AA)$, for $(e_1, e_2) \in EG^{[2]}$:
  \[
  \begin{array}{lll}
  &&(\pi_2^*-\pi^*_1)CS_\phi (\AA) (e_1, e_2)\\[2mm]
 &=& CS_\phi (\AA) (e_2) -CS_\phi (\AA) (e_1)\\[2mm]
 &=& CS_\phi (\AA) (e_1\cdot \hat{g}(e_1, e_2))-CS_\phi (\AA) (e_1)\\[2mm]
 &=& \hat{g}^* \omega_\phi - d\Phi (\AA, d\hat{g}\cdot \hat{g}^{-1}),
 \end{array}
 \]
 from which we can see that $\omega_\phi$ is a closed 3-form $\omega_\phi$ on $G$
 with periods in $\ZZ$. Hence, we can choose a bundle gerbe connection and curving on the multiplicative bundle
 gerbe $\cG$ such that the bundle gerbe curvature is given by
 $\omega_\phi$. Moreover, we can choose a bundle gerbe connection and curving
 on $(\cQ, X, EG^{[2]})$ whose bundle gerbe curvature is given by
 $(\pi_2^*-\pi^*_1)CS_\phi (\AA)$, hence the bundle 2-gerbe curving is given by
 the universal Chern-Simons form $CS_\phi (\AA)$.

 On the Deligne cohomology level, we obtain a degree 2 Deligne class
 in $H^2(G, \cD^2)$  associated to the multiplicative bundle gerbe $\cG$
 with connection and curving over $G$
and  to the degree 3 Deligne class in
 $H^3(BG, \cD^3)$ determined by $(\phi, \Phi (\disp{\frac{i}{2\pi}}F_\AA) )$, the
 Deligne class for the bundle 2-gerbe $\cQ$ over $BG$ with connection and curving.
 This completes the proof of the existence of
a universal Chern-Simons bundle 2-gerbe
 associated to $\phi$ and $(EG, \AA)$.
  \end{proof}

  The upshot of all this is that
 the universal Chern-Simons bundle 2-gerbe over $BG$ gives a geometric
  realization of the correspondence between three dimensional Chern-Simons gauge theories and
  Wess-Zumino-Witten models associated to $G$,
from which the proof of Theorem \ref{CS:image}
 immediately  follows.

  Let $P \rightarrow M$ be  a principal $G$-bundle with connection $A$. Let $f: M \rightarrow
  BG$ be a classifying map for this bundle with connection. This
  means that $f^* (EG,BG) \cong (P,M)$ and there exists a connection
  $\AA$ on $EG$ such that $f^* \AA = A$.

  \begin{definition}\label{csb2g} For the Chern-Simons
gauge theory canonically defined by
  a class $\phi \in H^4(BG, \ZZ) $,
  the {\it Chern-Simons bundle 2-gerbe} $\cQ_\phi(P, A)$
   associated with the principal $G$-bundle $P$
  with connection  $A$  over $M$ is defined to be the pullback of the universal Chern-Simons bundle
  2-gerbe by the classifying map $f$ of $(P, A)$.
  \end{definition}

  For a principal $G$-bundle $P\to M$ with a connection $A$, the corresponding
  Chern-Simons form for $(P, A)$ corresponding to the class $\phi \in H^4(BG, \ZZ)$
  is $$ CS_\phi(A) = f^*CS_\phi (\AA) \in \Omega^3(P),$$ such that
  \[
  dCS_\phi (A) =\pi^* \Phi (\disp{\frac{i}{2\pi}} F_A) \in \Omega^4_{cl, 0}(M).
  \]
  Hence, the curvature of the Chern-Simons bundle 2-gerbe $\cQ_\phi(P, A)$
  associated to $(P, A)$ over $M$
  is given by $\Phi (\disp{\frac{i}{2\pi}} F_A)$ and its bundle 2-gerbe curving is given by the
  Chern-Simons form $CS_\phi(A)$.

\begin{remark} We can connect our approach to the familiar Chern-Simons
action in the physics literature when $G=SU(N)$.
 We have a Deligne class constructed
using the pullback construction from the universal
Chern-Simons gauge theory
  \[
  c_\phi (P, A) \in H^3(M, \cD^3) \sta{(hol, curv)}{\longrightarrow}
  \check{H}(M, U(1)).
  \]
  Therefore for any smooth map $\sigma$ from a closed 3-dimensional
  manifold $Y$ to $ M$, the holonomy of $c_\phi (P, A)$ associated to
$\sigma$ is given by
  \ba\label{CS:action}
  e^{2\pi i CS_\phi (\sigma; A)} = hol (c_\phi (P, A)) (\sigma).
  \na
  See also \cite{Fre} for
  similar constructions. Now the Chern-Simons functional
  $CS_\phi(\sigma, A)$ will be the familiar formula
when $\Phi=cw^{-1}(r(\phi))$ is
 the second Chern polynomial for $SU(N)$. Fix  a   trivialisation of
  $(\sigma^*P, \sigma^*A)$ over $Y$,
   we can write the level $k$ Chern-Simons functional  as
  \[
  CS(\sigma, A) = \disp{\frac{k}{8\pi^2}\int_Y} Tr( \sigma^*A\wedge  \sigma^*dA + \disp{\frac 13}  \sigma^*A\wedge
  \sigma^* A \wedge \sigma^*A).
  \]
\end{remark}

  \begin{theorem}\label{csb2g=cs-invariant}
   With the canonical isomorphism between the Deligne cohomology
  and Cheeger-Simons cohomology,  the Chern-Simons bundle 2-gerbe $\cQ_\phi(P, A)$
  is equivalent in   Deligne cohomology to the Cheeger-Simons invariant
  $S_{\Phi, \phi} (P, A)$ described in Remark \ref{S:P-A}.
  \end{theorem}

  \begin{proof}
  It is well known that Cheeger-Simons differential characters are classified
  by Deligne cohomology (see, for example, \cite{Bry}). Stable
  isomorphism classes of bundle 2-gerbes
  with connection and curving are also classified by
Deligne cohomology (Proposition
  \ref{b2g:class}). The theorem of Cheeger and Simons given above defines a unique
  differential character satisfying certain conditions, thus it uniquely
  defines a class in Deligne cohomology and an equivalence class of bundle
  2-gerbes with connection and curving, so we must show that the corresponding
  Deligne class associated to  the Chern-Simons
  bundle 2-gerbe satisfies the required conditions.
  \begin{enumerate}
  \item The image under the map
  $H^3(M,\cD^3) \rightarrow \Omega_0^4(M,\RR)$ is the curvature 4-form of
  the Chern-Simons bundle 2-gerbe. The curvature of the universal CS bundle 2-gerbe is
  given by $\Phi(\disp{\frac{i}{2\pi}}F_\AA)$.
  Under pullback this becomes $\Phi(\disp{\frac{i}{2\pi}}F_A)$.
  \item The image of the map $H^3(M,\cD^3)
  \rightarrow H^4(M,\ZZ)$ is the characteristic 4-class of a bundle 2-gerbe.
  For the Chern-Simons bundle 2-gerbe, this shall be the pull-back of the
  4-class of the universal Chern-Simons bundle 2-gerbe by the classifying map, which
  is by construction $\phi \in H^4(BG,\ZZ)$.
  \item If a principal $G$-bundle $P_1$ with connection $A_1$,  is related to another
  principal $G$-bundle
  with connection $(P_2, A_2)$, via a bundle morphism $\psi$ then
  the corresponding classifying maps are related by $f_{P_1, A_1} = \psi \circ
  f_{P_2, A_2}$, so using both sides to pull back the universal CS bundle 2-gerbe
  we see that their corresponding Deligne classes behave
  as their Cheeger-Simons invariants $S_{\Phi, \phi}(P_1, A_1) = \psi^* S_{\Phi,\phi}(P_2, A_2)$.
  \end{enumerate}
  \end{proof}
  
  Recall that for a connection
  on a line bundle over $M$, a gauge transformation is given
  by a smooth function $M \to U(1)$, and an extended 
  gauge transformation for a bundle gerbe with connection and curving
  is given by a line bundle with connection over $M$.  We can discuss
  extended gauge transformations for a bundle 2-gerbe with connection
  and bundle 2-gerbe curving over $M$.
  
  The Chern-Simons bundle 2-gerbe $\cQ_{\phi}(P, A)$ is equipped with
  a bundle 2-gerbe  connection and curving such that the 2-curving 
  is given by the Chern-Simons form $CS_\phi (A)$ with 
  $$
  dCS_\phi (A) =\pi^* \Phi (\disp{\frac{i}{2\pi}} F_A) \in \Omega^4_{cl, 0}(M).
  $$
  Choose a covering $\{U_i\}$ of $M$, such that over
  $U_i$,  $\pi: P\to M$ admits  a section 
  $s_i$, then we obtain a \v{C}ech representative of the
  Deligne class corresponding to $\cQ_{\phi}(P, A)$
  \[
 c_\phi (P, A) =  [(g_{ijkl}^{\ }, A_{ijk}, B_{ij}, C_i)]
  \]
  with $C_i = s_i^* (CS_\phi (A))$. These local 3-forms
  $\{C_i\}$ are called the `{\sl C-field}' in string theory. 
  We can say that our Chern-Simons 
  bundle 2-gerbe $\cQ_{\phi}(P, A)$ carries the Chern-Simons form
  as the `{\sl C-field}'. As the bundle 2-gerbe curving (`C-field') is
  not uniquely determined, different choices are related by an
  extended gauge transformation. Suppose that
  $(g_{ijkl}^{\ }, A_{ijk}, B_{ij}, C_i)$ represents the Deligne
  class of $\cQ_{\phi}(P, A)$, then adding 
  a term $(1, 0, 0, \omega)$ with a closed 3-form of integer period
  $\omega \in \Omega^3_{cl, 0}(M)$ doesn't change
  the Deligne class, following from the exact sequence (\ref{exact:1}). 
  If $M$ is 2-connected,  we know
  that $\omega \in \Omega^3_{cl, 0}(M)$ canonically defines a
  bundle gerbe with connection and curving over $M$ whose
  bundle gerbe curvature is given by $\omega$. We call 
  this bundle gerbe an extended gauge transformation of
  the Chern-Simons bundle 2-gerbe. 
  
  Note that given another connection $A'$ on $\pi: P\to M$, the
  Chern-Simons bundle 2-gerbe $\cQ_{\phi}(P, A')$ is stably isomorphic
  to  $\cQ_{\phi}(P, A)$ as bundle 2-gerbes over $M$ ( they have the same
  characteristic class determined by $\phi$).  $\cQ_{\phi}(P, A)$
  and $\cQ_{\phi}(P, A')$ have
  different bundle 2-gerbe curving, on the level of Deligne cohomology,
   the difference is given   by 
   \[
   c_\phi (P, A) - c_\phi(P, A') = [( 1, 0, 0, CS_\phi(A, A'))]
   \]
   where $CS_\phi(A, A')$ is the well-defined Chern-Simons 3-form on $M$
   associated to a straight line path of connections on $P$  connecting $A$ and $A'$, 
   (we remind that  $CS_\phi(A)$ is  well-defined only   on $P$ in general).  It is natural
   (Cf.\cite{DFM}) to define the so-called {\sl C-field} on $M$ to be a pair
   $(A, c)$  where $A$ is a connection on $P$ and $c\in \Omega^3(M)$. So
   the space of $C$-fields is 
   \[
   \cA_P \times \Omega^3(M)
   \]
   where $\cA_P$ is the space of connections on $\pi: P\to M$. A $C$-field
   $(A, c)$ canonically defines a degree 3 Deligne class 
   \[
   c_\phi(P, A) + [(1, 0, 0, c)] \in H^3(M, \cD^3)
   \]
   through the Deligne class $c_\phi (P, A)$ of the
   Chern-Simons bundle 2-gerbe $\cQ_\phi (P, A)$. 
   The gauge transformation group for the space of $C$-fields  is defined to be
   \[
   \Omega^1(M, ad P) \times H^2(M, \cD^2),
   \]
   with the action on $C$-fields given by
   \ba\label{gauge:C}
   (\alpha, \cD) \cdot (A, c) = 
   \bigl(A + \alpha, c + CS_\phi (A, A + \alpha) + curv (\cD)\bigr)
   \na
   where $(\alpha,  \cD)\in  \Omega^1(M, ad P)\times H^2(M, \cD^2)$,
   and $curv (\cD) \in \Omega^3_{cl, 0}(M)$ is 
   the curvature of the degree 2 Deligne class $\cD$
   (or the corresponding bundle gerbe with connection and curving) on $M$.
    Then it is easy to see that two
   $C$-fields that are gauge equivalent under (\ref{gauge:C})
     define the same degree 3 Deligne class on $M$
   through the Chern-Simons bundle 2-gerbe, hence,
   the same Cheeger-Simons differential character on $M$.
   
   \begin{remark} Denote by  $\cG (P)$ the gauge transformation group of
    $P$ which   acts on $\cA_P \times \Omega^3(M)$ via 
   $    g \cdot (A, c) = 
   \bigl(A^g, c + CS_\phi (A, A^g) \bigr).
   $
   Due to the fact that
   \[
   CS_\phi (A, A^{g_1g_2}) - CS_\phi (A, A^{g_2}) - CS_\phi (A^{g_2}, A^{g_1g_2})
   \]
     depends on the choice of $A$, this $\cG (P)$-action is not a group action. It is
     observed in \cite{DFM} that if one interprets the space of
     $C$-fields with gauge group action (\ref{gauge:C}) as an
     action groupoid, then $\cG (P)$-action  is a sub-groupoid action.
     \end{remark}

  \section{Multiplicative Wess-Zumino-Witten models}\label{multi:WZW}

  In this section, we study the Wess-Zumino-Witten models in the
  image of the correspondence from $CS(G)$ to $WZW(G)$.

  Let $\cG$ be a bundle gerbe with connection and curving
   over a compact Lie group $G$, whose Deligne class
is in $H^2(G, \cD^2)$. With the identification between the Deligne cohomology
  and the Cheeger-Simons cohomology (\ref{del=cs}),
  \[
  (hol, curv):  H^2(G, \cD^2) \to \check{H}^2(G, U(1))
  \]
  where $hol$ and $curv$ are the holonomy and the curvature maps for the Deligne cohomology,
  we will define the bundle gerbe holonomy for stable equivalence classes of bundle
  gerbes with connection and curving.

  Let $\sigma: \Sigma \to G$ be a smooth map from a closed 2-dimensional
  surface $\Sigma$ to $G$. $\sigma$ represents a smooth 2-cocycle in $Z_2(G, \ZZ)$.
  Define  the holonomy of the bundle gerbe $\cG$ over $G$ to be the
  holonomy of the corresponding  Deligne class in $H^2(G, \cD^2)$, denoted
  by $hol_{\cG}$, then
  \[
  hol_{\cG} (\sigma) \in U(1),
  \]
  is called the bundle gerbe holonomy $hol_\cG (\cdot)$,
  evaluated on $\sigma: \Sigma\to G$.

  We point out that as $H^2(G, \cD^2)$ classifies stable isomorphism classes
  of bundle gerbes over  $G$ with connection and curving our
  bundle gerbe holonomy $hol_{\cG} (\sigma) $ depends only on the stable isomorphism
  class of $\cG$.

  \begin{proposition}\label{multiplicative:gerbe}
   For a  multiplicative bundle gerbe $\cG$ with connection and curving
   over $G$,
  the bundle gerbe holonomy satisfies the following multiplicative property:
  \ba
  \label{multiplicative:hol}
  hol_\cG (\sigma_1 \cdot \sigma_2 ) = hol_\cG (\sigma_1) \cdot hol_\cG(\sigma_2 )
  \na
  for any pair of smooth maps  $(\sigma_1, \sigma_2)$ from  any closed surface
  $\Sigma$ to $G$.
   Here $\sigma_1 \cdot \sigma_2$ denotes the smooth map
  from $\Sigma$ to $G$ obtained from the pointwise multiplication
  of $\sigma_1$ and $\sigma_2$with respect to
  the group multiplication on $G$.
  \end{proposition}
  \begin{proof} Recall our correspondence map in
Definition \ref{CS2WZW:map} and our
  integration map (\ref{integration:1}). We constructed $\Psi$ from
  a canonical $G$-bundle over $S^1\times G$  in Definition \ref{CS2WZW:map}
  such that the bundle gerbe holonomy for  the multiplicative bundle gerbe $\cG$ with
  connection and curving over $G$,
  being in  the image of $\Psi$, corresponds to the bundle 2-gerbe holonomy
  for the Chern-Simons bundle 2-gerbe $\cQ$ over $S^1\times G$ as follows.
Given a smooth map 
  $\sigma: \Sigma \to G$, $Id\times \sigma$ defines a smooth map
  $S^1\times \Sigma \to S^1\times G$, and
  \ba\label{hol:cG}
  hol_{\cG}(\sigma) = Hol_{\cQ}(Id\times \sigma),
  \na
  where $Hol_{\cQ}$ denotes the bundle 2-gerbe holonomy for the Chern-Simons bundle
  2-gerbe over $S^1 \times G$.

  Given a pair of smooth maps $\sigma_1$ and $\sigma_2$ from any closed surface
  $\Sigma$ to $G$. Denote by $\Sigma_{0, 3}$ a fixed sphere with three holes.
  We can construct a flat $G$-bundle over
  $\Sigma_{0, 3} \times \Sigma$ with boundary orientation given in such a way that
  the usual holonomies for flat $G$-bundle are $\sigma_1$, $\sigma_2$ and
  $\sigma_1\cdot \sigma_2 $ respectively. The Chern-Simons bundle 2-gerbe
  associated to this flat $G$-bundle is a flat bundle 2-gerbe in the sense that
  the bundle 2-gerbe holonomy is a homotopy invariant
(as follows from the exact sequence
  (\ref{exact:2})). This implies that
  the  bundle 2-gerbe holonomies for
   this flat Chern-Simons bundle 2-gerbe satisfy
   \ba\label{hol:cQ}
   Hol_{\cQ}(Id\times \sigma_1\cdot \sigma_2)
   =Hol_{\cQ}(Id\times \sigma_1) \cdot Hol_{\cQ}(Id\times \sigma_2).
   \na
   Combining (\ref{hol:cG}) and (\ref{hol:cQ}), we obtain the multiplicative
   property for the  bundle gerbe holonomy of $\cG$:
   \[
   hol_\cG (\sigma_1 \cdot \sigma_2 ) = hol_\cG (\sigma_1) \cdot hol_\cG(\sigma_2 )
  \]
  for any pair of smooth maps  $(\sigma_1, \sigma_2)$ from  any closed surface
  $\Sigma$ to $G$.
  \end{proof}

  Recall that $d_i$ ($i=0, 1, 2$) are the face 
  maps from $G\times G \to G$ such that
  $d_0(g_1, g_2) = g_2$, $d_1(g_1, g_2) = g_1g_2$ and $d_2(g_1, g_2) = g_1$ for
  $(g_1, g_2) \in G$. Let $\cG$ be a bundle gerbe with connection and curving.
   Let $curv(\cG)$ be the bundle gerbe curvature of $\cG$.
  We can consider the pair $(\sigma_1, \sigma_2)$ as a map into $G \times G$. Hence
  we can define the holonomy of $(\sigma_1, \sigma_2)$ with respect to the
  bundle gerbe connection and curving on
  \[
  \delta(\cG) = d_0^*(\cG)\otimes d_1^*(\cG^*)\otimes d_2^*(\cG),
  \]
  whose curvature is given by
  \[
  d_0^*(curv(\cG)) - d_1^*(curv(\cG)) + d_2^*(curv(\cG)) = dB,
  \]
  for a 2-form $B$ on $G \times G$.

  Because we have a trivialisation for $\delta(\cG)$,
  we can calculate this holonomy as
  $$
  \exp \int_{\Sigma} (\sigma_1, \sigma_2)^*B .
  $$
  On the other hand we can compose $(\sigma_1, \sigma_2)$
  with the three maps into $G$. This gives $\sigma_1$, $\sigma_1\cdot\sigma_2$
  and $\sigma_2$ where the second of these maps is the result of pointwise
  multiplying in $G$.  Because
  the bundle gerbe connection and curving on
  \[
  \delta(\cG) = d_0^*(\cG)\otimes d_1^*(\cG^*)\otimes d_2^*(\cG)
  \]
   are $\delta$ of those on $\cG$,  then the bundle gerbe  holonomy can also be calculated as
  $$
  hol_\cG(\sigma_1) hol_\cG(\sigma_1\sigma_2)^{-1} hol_\cG(\sigma_2)
  $$
  so that we have
  $$
  hol_\cG(\sigma_1)  hol_\cG(\sigma_2) = \left(\exp \int_{\Sigma} (\sigma_1,
  \sigma_2)^*B\right) hol_\cG(\sigma_1\cdot \sigma_2) .
  $$
The following proposition gives a
  necessary condition for a bundle gerbe with connection and curving
  to be multiplicative, and can be proved by direct calculation.

  \begin{proposition}
  \label{criteria} If a  bundle gerbe $\cG$ with connection
  and curving on $G$ is multiplicative, then
  $d_0^*\cG \otimes d_2^*\cG$ is stably isomorphic to $d_1^* \cG$ as bundle gerbes
  over $G\times G$, and there exist a imaginary valued 2-form $B$ on $G\times G$
  such that
  \[\begin{array}{c}
  d_0^*(curv(\cG)) - d_1^*(curv(\cG)) + d_2^*(curv(\cG)) = dB, \\[2mm]
  \disp{\int_{\Sigma}} (\sigma_1, \sigma_2)^*B \in 2\pi i \ZZ,\end{array}
  \]
  for any pair of smooth maps  from any closed surface
  $\Sigma$ to $G$.
  \end{proposition}

  For a  Wess-Zumino-Witten model with group manifold $G$ in the image of the
  correspondence map $\Psi: CS(G)\to WZW(G)$, we know that the bundle gerbe
  $\cG$ with connection and curving on $G$ is multiplicative.
The Wess-Zumino-Witten action
regarded as a function on
the space of smooth maps $\{\sigma: \Sigma\to G \}$ exponentiates to
the bundle gerbe holonomy of $\cG$, that is, for a smooth map $\sigma$, we have
  \[
 exp\bigl( S_{\text{wzw}} (\sigma) \bigr)= hol_{\cG} (\sigma).
  \]
  From Proposition \ref{multiplicative:gerbe}, we know that the Wess-Zumino-Witten
  action for a multiplicative Wess-Zumino-Witten model
  satisfies the following property:
  \[
 exp \bigl(S_{\text{wzw}}(\sigma_1\cdot \sigma_2) \bigr)
 = exp \bigl(S_{\text{wzw}} (\sigma_1)\bigr)\cdot
 exp \bigl(S_{\text{wzw}} (\sigma_2)\bigr),
  \]
  for  a pair of smooth maps $\sigma_1$ and $\sigma_2$ from any closed surface
  $\Sigma$ to $G$.

  From the commutative diagram:
  \[
  \begin{array}{ccc}
  H^3(BG, \cD^3)&\longrightarrow& H^2(G, \cD^2)\\[2mm]
\downarrow^c && \downarrow^c \\[2mm]
H^4(BG, \ZZ)& \to & H^3(G, \ZZ)
\end{array}
\]
 we see that for a general compact semi-simple Lie group $G$,
 $H^3(BG, \cD^3)\to  H^2(G, \cD^2)$ is not surjective. In particular,
 for non-simply connected compact semi-simple Lie group $G$, we
 know that the Wess-Zumino-Witten model on $G$ is
only multiplicative at certain levels. For example, the 
Wess-Zumino-Witten model on $SO(3)$ is multiplicative if and   only if
the Dixmier-Douady class of the corresponding bundle gerbe is
an even class in $H^3(SO(3), \ZZ) \cong \ZZ$.

  \end{document}